\documentstyle[12pt,amssymb]{article}

\newcommand{\no}{\noindent}
\def\C{{\Bbb C}}
\def\H{{\Bbb H}}
\def\N{{\Bbb N}}
\def\Q{{\Bbb Q}}
\def\R{{\Bbb R}}
\def\P{{\Bbb P}}
\def\Z{{\Bbb Z}}
\def\CC{{\cal C}}
\def\FF{{\cal F}}

\def\JJ{{\cal J}}
\def\YY{{\cal Y}}
\def\VV{{\cal V}}
\def\XX{{\cal X}}

\def\extra{\def\normalbaselines{\baselineskip 20pt\lineskip
3pt\lineskiplimit 3pt}} \def\ex{\def\normalbaselines{\baselineskip
18pt\lineskip 3pt\lineskiplimit 3pt}}
 \def\1{{\rm l}\hskip
-0.21truecm 1}

\newcommand{\bi}{\begin{itemize}}
\newcommand{\ei}{\end{itemize}}
\newcommand{\be}{\begin{enumerate}}
\newcommand{\ee}{\end{enumerate}}
\newcommand{\bc}{\begin{center}}
\newcommand{\ec}{\end{center}}
\newcommand{\ba}{\begin{array}}
\newcommand{\ea}{\end{array}}
\newcommand{\bq}{\begin{quote}\small}
\newcommand{\eq}{\end{quote}}

\begin{document}

\bc
{\LARGE Hodge numbers attached to a\\
polynomial map}\\
\bigskip

by\\
\medskip
{\sl R. Garc\'\i a L\'opez\quad{\rm and}\quad A. N\'emethi}\\
\ec
\vspace*{1cm}

\no {\bf 1. Introduction}
\bigskip

\no (1.1) In the last years the behaviour at infinity of polynomial maps
has been extensively studied, basically by means of topological
invariants: Eisenbud-Neumann diagrams, monodromy at infinity,~etc.
In order to carry further this analysis of the asymptotic behaviour of
polynomial maps we consider in this paper the following approach: One
can attach to any polynomial map a geometrical variation of mixed Hodge
structures (from now on, MHS; see (1.2) below for details).  Using
either the geometric methods of \cite{St-Zu}, \cite{Na} or Saito's
theory of mixed Hodge modules (\cite{Sa}), one gets that this variation
gives rise to a functorial limit MHS. The equivariant Hodge numbers of
this limit MHS are analytical invariants of the polynomial map under
consideration and in this paper we compute them for a class of
polynomials which was extensively studied in \cite{GN1}, \cite{GN2}
from a topological point of view.  More precisely, we determine the
equivariant Hodge numbers of this limit MHS in terms of
equivariant Hodge numbers attached to some isolated
hypersurface singularities and Hodge numbers of cyclic coverings
of projective space branched along a hypersurface (see theorem (3.1)
for the precise statement). Both these local
and global Hodge numbers can be explicitly computed in a number of
cases, as shown in the examples in section 7. \medskip

\no (1.2) We give the precise description of the invariants we are
going to study. Let $f:\C^{n+1}\to\C$ be a polynomial map.  It is known
that there is a finite set $\Gamma_f$ such that $f$ defines a locally
trivial $\CC^\infty$-fibration over $\C-\Gamma_f$ (cf. e.g.
\cite[Appendix]{Ph}). Given a positive
number $r\in\R$, denote by $D_r$ the disk in the complex plane of center
$0$ and radius $r$, and set $D^*_r=D_r-\{0\},
Z_r=\C^{n+1}-f^{-1}(D_r)$.  If $r>\max\{\,\vert x\vert\ :\
x\in\Gamma_f\,\}$ then the map
$$
1/f:Z_r\longrightarrow
D^*_{1/r}
$$
is a locally trivial $\CC^{\infty}$-fibration.  Actually we
have a projective system of equivalent fibrations indexed by $r$, thus
we can assume $r$ so big as necessary.

\no Now we compactify the map $1/f$ and we add a fiber over $0\in\C$.
That is, we take an analytical manifold $X$ where $Z=Z_r$ is embedded as an
open dense subset, $(X-Z)_{red}$ is a divisor with normal crossings and
smooth
irreducible components and the map $1/f$ extends to a flat projective
map $p:X\longrightarrow D_{1/r}$, smooth over $D^*_{1/r}$.  Set
$X-Z=Y\cup\Delta$, where $Y=p^{-1}(0)$ and $\Delta$ is the union of the
irreducible components of $X-Z$ not in~$Y$.  Increasing $r$ if
necessary, we can assume that the restriction of $p$ to any intersection
of components of $\Delta _{red}$ is smooth, so that $\Delta$ becomes a
relative divisor with normal crossings.
Let $\H$ be the universal covering space of $D^*_{1/r}$,
set $\tilde
Z=Z\times_{D_{1/r}}\H$.  Notice that $\tilde Z$ has the homotopy type of
a generic fiber of~$f$.
Then
the cohomology groups $H^i(\tilde Z,\Q)$ carry a limit MHS. This follows
for example from \cite[\S 5]{St-Zu} or from
\cite[Th\'eor\`eme~5.13]{Na}. All MHS we will work with will be
regarded as $\Q$-MHS.\\
One can prove using standard arguments that this MHS does not
depend on the chosen compactification and is functorial
on $f$ (and therefore it is for example an invariant of polynomials with
respect to algebraic changes of variables in $\C^{n+1}$).
\medskip

\no (1.3) {\it Definition}: Given a polynomial map
$f:\C^{n+1}\longrightarrow\C$, the limit MHS on $H^n(\tilde Z,\Q)$
obtained in the way described above will be called the MHS at infinity
of~$f$. \medskip

\no (1.4) \ The MHS at infinity
should be viewed as a global analogue of the limit MHS associated to a
hypersurface singularity in \cite{St.Oslo} and \cite{Na}.  If $f$ is a
cohomologically tame polynomial (see \cite{Sab}) there is another
definition of a MHS at infinity due to C.~Sabbah.
For these polynomials, Sabbah has proved that
the equivariant Hodge
numbers of the two MHS's are the same.
Unless otherwise stated, all (co)homology
groups appearing in this paper will be assumed
to be reduced and to have complex coefficients.
We choose once and for all a root $i$ of $-1$ in $\C$ (equivalently, an
orientation
on the complex plane) and from now on we will use the notation
$e(r):=e^{2 \pi i r},\ r\in \R$.

\bigskip

\no {\bf 2. The main construction and exact sequences}
\bigskip

\no (2.1) Although the MHS at infinity is defined for any
polynomial map, if one does not place some kind of restriction on $f$
even the topological behaviour at infinity can be very complicated.
\smallskip

\no We will consider the class of polynomials
$f:\C^{n+1}\longrightarrow\C$ which satisfy the following condition:
\goodbreak
$$
(*)\cases{\hbox{If
$t\in\C$ is not a critical value of $f$,}\cr \hbox{then the closure in
$\P^{n+1}$ of the affine}\cr \hbox{hypersurface $\{f=t\}$ is
non-singular.}}
$$
Polynomials satisfying this condition will be called
$(*)$-polynomials.  In this section we recall some consequences of the
$(*)$ condition and some exact sequences attached to a $(*)$-polynomial.
For details, see \cite{GN1}, \cite{GN2}.
Let $f:\C^{n+1}\longrightarrow\C$ be a $(*)$-polynomial. Then:\medskip

\no (a) Any fiber of $f$ has at most isolated singularities and has the
homotopy type of a bouquet of $n$-dimensional spheres.  Consider the
fibration at infinity of $f$:
$$
f:f^{-1}(\partial\,\bar D_r)\cap
B_R\longrightarrow\partial\,\bar D_r\, ,
$$
where $D_r$ (resp. $B_R$) is a
disc in $\C$ (resp. a ball in $\C^{n+1}$) of radius $r$ (resp. $R$) and
$0<<r<<R$.
This fibration is equivalent to the fibration
$f:f^{-1}(\partial\,\bar D_r)\longrightarrow \partial\,\bar D_r$.  For
$t\in\partial\,\bar D_r$ set $X^0_t=f^{-1}(t)\subseteq\C^{n+1}$.
The geometric monodromy of this fibration (where  $\partial\,\bar D_r$
is positively oriented) is denoted by $T_{geom}$.
In \cite{GN1,GN2} we proved some  basic properties of the {\em classical}
algebraic monodromy
$T^\infty_f:H^n(X^0_t)\longrightarrow H^n(X^0_t)$ which is defined by
$T^\infty_f[\omega]=[T_{geom}^*(\omega)]$.

\no Let $X_t$ denote the projective closure of the affine
hypersurface $X^0_t$. The {\em classical} algebraic
monodromy over $\partial\,\bar D_r$ on the $n$-th cohomology groups of
the projective
closures will be denoted by $T:H^n(X_t)\longrightarrow H^n(X_t)$
(and is defined also as a pull--back).
\smallskip

\no (b) For $R>>0$, the Milnor fibration at infinity
$\varphi=\frac{f}{\Vert
f\Vert}:\partial\,\bar B_R-f^{-1}(0)\longrightarrow S^1$ exists and it
is equivalent to the fibration at infinity of $f$.  In particular,
there is a well defined variation map $Var: H^n(X^0_t)\longrightarrow
H_c^n(X_t^0)$ ($Var([\omega]):=[T_{geom}^*(\omega)-\omega]$),
which is an isomorphism and has the usual compatibility
properties with the monodromy $T^\infty_f$. \smallskip

\no (c) The projective hypersurface $X^\infty:=\{f_d=0\}\subseteq\P^n$
has
only isolated singularities.  The semisimple part of $T^\infty_f$
depends only on the local topological type of these singular points, but
the nilpotent part of $T^\infty_f$ (and also the intersection form on
$H_n(X_t^0)$) depends also on the position of the singularities of
$X^\infty$. \smallskip

\no In general, if one wants to relate the algebraic monodromy with
Hodge
theory, it is more convenient to use the algebraic monodromy induced on
the (horizontal sections of the) cohomological Gauss--Manin connection
(see e.g. \cite{AGV} \S 12).  This is the inverse dual of the classical
homological monodromy [loc. cit 13.1.A].  In this paper they will be
denoted by $M^{\infty}_f=(T^{\infty}_f)^{-1}$ acting on $H^n(X^0_t)$,
$M=T^{-1}$ acting on $H^n(X_t)$; and we will call them simply
(algebraic) monodromies.  Obviously, $W([\omega]):=
[(T_{geom}^*)^{-1}(\omega)-\omega]$ has the well--known compatibility
with $M^{\infty}_f$ (see 2.3).  We will use the corresponding similar
notations in all local situations which will appear in the body of the
paper.

\smallskip

\no (2.2) In order to study the MHS at infinity of a $(*)$-polynomial we
will consider a concrete compactification of the map $1/f$ (see
\cite{GN1}, \cite{GN2}). \medskip

\no Let $f=f_d+f_{d-1}+\dots$\ \ denote the decomposition of $f$ into
homogeneous components, let $D$ denote a disc in the complex plane with
center at the origin and sufficiently small radius.  Set
$$
\XX=\{\,([z],
t)\in\P\,^{n+1}\times D\ \vert\ t(f_d+z_0\,f_{d-1}+\dots +z_0^d\,
f_0)=z^d_0\,\}
$$
where $[z_0:z_1:\dots :z_{n+1}]$ are homogeneous coordinates in
$\P\,^{n+1}$.

\no The map $\pi:\XX\longrightarrow D$ given by $\pi([z], t)=t$ induces
a
locally trivial $\CC^\infty$ filtration over $D-\{0\}$ and the fibers of
$\pi$ are exactly the projective closures of the fibers of~$f$.  The
monodromy of $\pi$ over $\partial D$ (with positive orientation) is
$M$.

\no Then by \cite{GN1} we have a commutative diagram
$$
\ex\matrix{&0&\to&P^n_{X^\infty}(X_t)&\to&P^n(X_t)&\to&H^n(X^0_t)
&\to&P_{X^\infty}^{n+1}(X_t)&\to&0\cr \hbox{\small\rm (E.0)}
&&&\downarrow id&&\downarrow
M&&\downarrow M_f^\infty&&\downarrow id\cr
&0&\to&P^n_{X^\infty}(X_t)&\to&P^n(X_t)&\to&H^n(X^0_t)&\to&P_{X^\infty}
^{n+1}(X_t)&\to&0}
$$
where $P$ denotes primitive cohomology
(i.e. for $X\subset \P\,^{N}$ one has
$P^k(X)=$
{\it Coker}$(H^k(\P^N)\to H^k(X))$).  The above spaces carry natural
MHS's: $P^n(X_t)$ carries Schmid's $MHS,\ H^n(X^0_t)$ the MHS introduced
in section 1 an all arrows are morphisms of MHS's. \medskip

\no The compactification provided by $\pi$ is not such a good one since
its center fiber $\pi^{-1}(0)$ is not reduced (it has multiplicity $d$)
and $\XX$ is very singular.
Therefore, we will consider the normalization of the $d$-fold covering
of~$\pi$.  More precisely, let $D'$ be again a small disk,
$\delta:D'\longrightarrow D$ the map given by $\delta(t)=t^d$. Define
${\XX}'$ as the normalization of $\XX\times_\delta D'$ and let
$\pi':{\XX}'\longrightarrow D'$ be the natural projection.  The map
$\pi'$
gives a fibration over $D'-\{0\}$ with the same fiber as $\pi$ and
classical monodromy $T^{-d}$.  Both $\XX'$ and the central fiber
$X_0'=(\pi')^{-1}(0)$ have only isolated singularities.  Set~${\rm
Sing}\,X^\infty=\{p_1,\dots ,\,p_k\}$ and let $F_j, M_j$ be the local
Milnor fiber and the algebraic monodromy
$M_j:H^{n-1}(F_j)\longrightarrow H^{n-1}(F_j)$ of the hypersurface
singularity $(X^\infty, p_j)\subset (\P\,^n,p_j)$ for $1\le j\le k$
(Here $\P\,^n$ is the hyperplane at infinity). \medskip

\no Then ${\rm Sing}\,{\XX}'= {\rm Sing}\, X^\infty\times \{0\}$
and the central
fiber $X'_0$ is the $d$-fold cyclic covering of $\P^n$ branched along
$X^\infty$.  In particular, if we set ${\rm
Sing}\,X'_0=\{p'_1,\dots,\,p'_k\}$ then the isolated singularities
$(X'_0, p'_j)$ are the $d$-th suspensions of the singularities
$(X^\infty, p_j)$.

\no Let $F'_j$ (resp. $M'_j$) be the Milnor fiber (resp. the monodromy)
of
the local smoothing $\pi':(\XX', p'_j)\longrightarrow (D', 0)$. We
remark that this smoothing is
not the (natural) smoothing of $(X'_0, p'_j)$ with total space
smooth. The following
commutative diagram is the exact sequence of vanishing
cycles corresponding
to $\pi'$ together with the monodromy action. It is a diagram of MHS's
(cf. \cite[(10)]{GN1} or \cite[(2.2)]{GN2}).

$$\ex\matrix{&0\!&\to&P^n(X'_0)&\to&P^n(X_t)&\to&
\displaystyle\!\!\!\bigoplus^k_{j=1}H^n\,(F'_j)\!\!&\to
&P^{n+1}(X'_0)&\to&\!0\cr
\hbox{\small\rm (E.1)}&&&\downarrow id&&\downarrow
M^{d}&&\downarrow\oplus\,(M'_j)^{-1}&&\downarrow id\cr
&0\!&\to&P^n(X'_0)&\to&P^n(X_t)&\to&
\displaystyle\!\!\!\bigoplus^k_{j=1}\,H^n(F'_j)\!\!&\to
&P^{n+1}(X'_0)&\to&\!0.}
$$

\no Notice that $X'_0\subset\P^{n+1}$ is given by the equation
$z_0^d=f_d(z_1,\dots,\,z_{n+1})$.  There is a Galois action of $\Z/d\Z$
on $X'_0$ generated by the automorphism $[z_0:\dots
:\,z_{n+1}]\longmapsto [e(1/d)\,z_0:\dots :\,z_{n+1}]$ which
induces at the cohomology level an automorphism which will be denoted
$G^q:P^q(X'_0)\longrightarrow P^q(X'_0)$.

\no The disadvantange of (E.1) is
that the monodromy action at the level of $P^n(X_t)$ is $M^{d}$ instead
of $M$. In \cite{GN2} we constructed an automorphism of $\XX'$ over $D'$
which
is a ``lifting" of the geometric monodromy of~$\pi$.  In [loc. cit.] we
identified $(H^n\,(F'_j),\,M'_j)$ with
$$
\Big(H^{n-1}\,(F_j)^{\oplus\,(d-1)},\ c_{d-1}\,(M_j)^d\Big)
$$
where for an operator $\varphi:V\longrightarrow V$ and an integer
$\ell\ge 1$
we denote $c_\ell(\varphi):V^{\oplus\ell}\longrightarrow V^{\oplus\ell}$
the operator defined by $c_\ell(\varphi)\
(x_1,\dots,\,x_\ell)=(\varphi(x_\ell),\break x_1,\dots,\ x_{\ell-1})$.
Now the ``lifted monodromy action" induces the following commutative
diagram (of vector spaces):
$$
\ex\matrix{&0\!\!\!&\to&P^n(X'_0)&\to&P^n(X_t)&\to&\displaystyle
\!\!\!\bigoplus^k_{j=1}H^{n-1}\,
(F_j)^{\oplus\,(d-1)}\!\!&\to&P^{n+1}(X'_0)&\to&\!\!\!0\cr
\hbox{\small\rm (E.1')}&&&\downarrow G^n&&\downarrow
M&&\downarrow\oplus\,c_{d-1}\,(M_j^{-1})&&\downarrow G^{n+1}\cr
&0\!\!\!&\to&P^n(X'_0)&\to&P^n(X_t)&\to&\displaystyle\!\!\!\bigoplus
^k_{j=1}\,
H^{n-1}(F_j)^{\oplus(d-1)}\!\!&\to&P^{n+1}(X'_0)&\to&\!\!\!0.}
$$
(In [loc. cit.] it is not clarified the MHS-meaning of this diagram,
this will be done in the next section.  That discussion will show that
(E.1') is actually a diagram of MHS's. One should be aware of the
fact that if we consider
 on
$H^n(F_j')$ and $H^{n-1}(F_j)$ their natural MHS, then
$H^n(F_j')\not =H^{n-1}(F_j)^{\oplus(d-1)}$ as MHS's.)

\no Theorem 2 in the Appendix of \cite{GN1} together with
 the ``lifted monodromy action"
gives the following commutative diagram of MHS (cf. \cite{GN2}):
$$
\ex\matrix{&0\!\!\!&\to&P^n(X'_0)&\to&{\rm Ker}\,(M^{d}-
Id\ \vert\ P^n(X_t)\,)&\to&\displaystyle
\!\!\!\bigoplus^k_{j=1}H^{n+1}_{\{p'_j\}}\,(\XX')&\to&\!\!\!0\cr
\hbox{\small\rm (E.2)}&&&\downarrow G^n&&\downarrow
M&&\downarrow id\cr
&0\!\!\!&\to&P^n(X'_0)&\to&{\rm Ker}\,(M^{d}-
id\ \vert\ P^n(X_t)\,)&\to&\displaystyle
\!\!\!\bigoplus^k_{j=1}H^{n+1}_{\{p'_j\}}\,(\XX')&\to&\!\!\!0.}
$$
(2.3)\
There is a nice connection between the variation map and the
diagram \ (E.0)\  (for details, see \cite[pp.~217--218]{GN1}).
There is an orthogonal decomposition
$P^n(X_t)_1=P^n_{X^\infty}\,(X_t)\oplus P^n_{X^\infty}\,(X_t)^\bot,\
M_1=Id\oplus M'$ (where $P^n(X_t)_1$ denotes the space of generalized
eigenvectors of $P^n(X_t)$ of eigenvalue~1 with respect to $M$, we
will use the same notation for the other groups as well).

\no We identify $P^n_{X^\infty} (X_t)^\bot$ with the image of
$P^n\,(X_t)\longrightarrow H^n\,(X^0_t)$.  Then the composed map
$$
w:H^n\,(X^0_t)_1\buildrel {\rm W}_1\over\longrightarrow
H^n_c\,(X_t^0)_1\longrightarrow H^n\,(X^0_t)_1
$$
(where $W$ was defined in (2.1)) has range exactly
$P^n_{X^\infty}(X_t)^\bot$ and makes the following diagram (of vector
spaces) commutative
$$
\begin{array}{ccccccccc}
0\!\!\!&\to&P^n_{X^\infty}(X_t)^\bot&\to&\!\!\!H^n(X^0_t)_1
&\to&P^{n+1}_{X^\infty}
(X_t)&\to&\!\!\!0\\
&&&&\!\!\!\!\!\!\!\!\!\!\!\!\!\!\!\!\!\!w\!\! \swarrow\qquad
\qquad\qquad\\[-3.5mm]
&&\downarrow M'-id&&\!\!\!\!\!\!\swarrow\qquad\downarrow
(M^\infty_f)_1-id&&\downarrow 0\\[3mm]
0\!\!\!&\to&P^n_{X^\infty}(X_t)^\bot&\to&\!\!\!H^n(X^0_t)_1&
\to&P^{n+1}_{X^\infty}
(X_t)&\to&\!\!\!0.
\end{array}
$$
\smallskip

\no {\bf 3. The mixed Hodge structure at infinity in the $(*)$-case.}
\bigskip

\no In this section we describe the equivariant Hodge numbers of the MHS
at infinity of a $(*)$-polynomial.

\no Let $f$ be a $(*)$-polynomial and keep the notations of section 2.
Let $(M^\infty_f)_{ss}$ (resp. $(M^\infty_f)_u$) denote the semisimple
(resp. the unipotent) part of the monodromy at infinity $M^{\infty}_f$.
We recall (\cite[Th\'eor\`eme~15]{Na}) that $(M^\infty_f)_{ss}$ is a
MHS-automorphism of $H^n(X^0_t)$ of type $(0,0)$ and
$N^0=\log\,(M^\infty_f)_u$ is a MHS-automorphism of type $(-1,-1)$.
We will denote $H^n(X_t^0)_{\xi}$ the $\xi$-eigenspace of $H^n(X_t^0)$
with respect to $M^{\infty}_f$.
For $\xi$ a $d$-th root of unity, we denote by $P^\ell(X'_0)_\xi$ be
the $\xi$-eigenspace of $P^\ell(X'_0)$ with respect to the Galois action
$G^\ell$.  Obviously, this eigenspace decomposition is compatible with
the MHS of $P^\ell(X'_0)$.
\smallskip

\no We recall the definition of primitive spaces with respect to a
nilpotent endomorphism. (cf. e.g.
\cite[(5.3)]{Ne1}, \cite[(6.4)]{Sch}): Let $H$ be a finitely dimensional
$\Q$-vector space and $N:H\to H$ a nilpotent endomorphism of $H$. Let
$W$
denote the weight filtration on $H$ corresponding to $N$ and centered at
$m$. Assume there is a fitration $F$ on the complexification of $H$
such that $(H,W,F)$ is a MHS. For $r\ge 0$ one sets:
\[
P_r := \mbox{Ker} N^{r+1}: \mbox{Gr}^W_{m+r}H\to\mbox{Gr}^W_{m-r-2}H.
\]
Then $P^r$ carries a MHS of weight $m+r$. For $p+q\ge m$ we will call
$P^{p,q}_{p+q-m}$ the primitive space of type $(p,q)$.
\smallskip

\no The local vanishing cohomology $H^{n-1}(F_j)$ of $(X^\infty, p_j)$
has a
MHS (cf.~\cite{St.Oslo}).  For $\xi\ne 1$ (resp. for $\xi=1$) the
weight
filtration of the generalized $\xi$-eigenspace is the monodromy weight
filtration centered at $n-1$ (resp. at $n$).
In this case, we will  denote the dimension of the corresponding
primitive spaces (with respect to the infinitesimal nilpotent monodromy)
by $p^{pq}_\xi(X^\infty ,p_j)$.
The main theorem of this section is the following:
\medskip

\no (3.1) {\bf Theorem:}{\it
\smallskip

\noindent a) The weight filtration of $H^n(X^0_t)_1$ is the monodromy
weight filtration given by $N_1^0$ with center $n+1$.
The dimensions of the primitive spaces of type $(p,q)$
$(p+q\geq n+1)$ are:
$$
p^{p,q}_1(f)=h^{n-q,\ n-p}\ (P^{n-1}(X^\infty)\,).
$$
\noindent b)
The weight filtration of $H^n(X^0_t)_{\ne 1}$ is the monodromy weight
filtration given by $N^0_{\ne 1}$ with center $n$.  The dimensions of the
primitive spaces of type $(p,q)$ are given as follows: \smallskip

i) If $\xi^d=1$ and $\xi\ne 1$, then
$$p^{p,q}_{\xi} (f)=h^{n-q,\ n-p}\,(P^n(X'_0)_{\xi}).
$$
\smallskip

ii) If $\xi^d\ne 1$, then
$$
p^{p,q}_{\xi}(f)=\sum_{j=1}^k\
p^{p-[\beta+\gamma],\
q-\delta+[\beta+\gamma]}_{\xi^{1-d}}\ (X^\infty ,p_j)
$$
where
$\xi=e^{-2\pi i\beta}(0<\beta<1),\ \gamma=\{(d-1)\,\beta\},\
$ and $\delta=0$ if $\gamma=0$ and $\delta=1$ if
$\gamma\ne 0$.

\no $($ Here $[\quad]$ denotes integer part and for $\alpha\in\R
,\ \{\alpha\}=\alpha-[\alpha]\,)$. }  \medskip

\no Notice that for $\xi^d\ne 1,\ p^{p,q}_\xi(f)$ depends only on
the
type of the local singularities $(X^\infty, p_j)^k_{j=1}$ (and not on
their position). \medskip

\noindent (3.2){\it Remarks:\ }
a) The equivariant Hodge numbers $h^{a,b}_\lambda$ of
$H^n(X^0_t)$ can be computed from the formula
$$
h^{a,b}_\lambda=\sum_{\ell\ge 0}\ p^{a+\ell,\ b+\ell}_\lambda\
(f)
$$
where $a+b\ge$ center of the weight filtration $:=c$, and from
$h_\lambda^{a,b}=h_{\bar\lambda}^{c-a,\ c-b}
=h_{\lambda}^{c-b,\ c-a}$  for $a+b\leq c$ (cf. \cite[(5.3)]{Ne1}).
\smallskip

\noindent b) The above theorem shows that the weight filtration of
$H^n(X^0_t)$ has a behaviour similar to that of the weight filtration of
the MHS on the vanishing cohomology corresponding to an isolated
hypersurface singularity (namely, it is centered at $n$ for
$\lambda\neq 1$ and at $n+1$ for $\lambda=1$).
For $\lambda\ne 1$ this can be proved using Schmid's results \cite{Sch}
and for $\lambda=1$ follows basically from the fact that the variation
map is an isomorphism. \smallskip

\noindent c) It follows also from the theorem that the maximal size of a
Jordan block of $(M^\infty_f)_\xi$ is
$n+1$ if $\xi^d=1$ and $\xi\ne 1$, otherwise is $n$
(cf.~\cite[p. 4]{GN2}). \smallskip

\noindent d) If $f$ and $f'$ are two $(*)$-polynomials with the same
highest degree form, their equivariant Hodge numbers are the same.
\bigskip

\no {\sl Proof of theorem 3.1:}

\no We will say that two MHS's $H_1$ and $H_2$ are numerically
equivalent (and we will denote this $H_1\buildrel n\over=
H_2$) if their Hodge numbers are the same. The primitive decomposition
 theorem
gives the following fact, which will be used several times.
\medskip

\noindent (3.3) {\bf Lemma:}\
{\em Let $H_i$ ($i=1,2$) be a MHS with a nilpotent
$(-1,-1)$--morphism $N_i$ such that the weight filtration of $H_i$ is
the monodromy weight filtration of $N_i$ (with the same fixed center).
If ${\rm Ker} N_1\buildrel n\over=
{\rm Ker} N_2$, then $H_1\buildrel n\over=
H_2$. }
\smallskip

We divide the proof of (3.1) in the following cases.

\smallskip

\noindent {\bf Case $\xi=1$.}
\smallskip

\no In order to prove this case we study first the MHS's of
$P^*(X^\infty)$ and $P^*_{X^\infty}(X_t)$. (Some of the results of this
subsection will be used in section~4 too.) \medskip

\noindent (3.4) {\bf Lemma:}
$P^{n+i}_{X^\infty}(X_t)$ and $P^{n-i}(X^\infty)$ are dual MHS's
 with respect to $\Q(-n).$ $P^n(X^\infty)$ is pure of weight $n$ and
$P^{n-1}(X^\infty)=W_{n-1}\ P^{n-1}(X^\infty)$.
\medskip

\noindent{\sl Proof of the lemma}: The duality statement follows quite
directly from
Saito's formalism: On the category of mixed Hodge modules one can define
a duality functor which conmutes with the nearby cycle functor up to a
Tate twist (cf.~\cite{Sa}, 2.6) and extends to a functor defined in the
derived category, compatible with Verdier duality.
In Saito's terminology and with the notations of (1.2), if
$i:\Delta\hookrightarrow
X$ is the inclusion map one has to apply this duality to the
object $\R\Psi_p\,
i_*\,i^!\,\Q_X [n+1]$ in the derived category of mixed Hodge modules.
The purity of $P^n(X^\infty)$ is proved in~\cite[(4.7)]{St.Oslo}.
\smallskip

\no (3.5) Consider now a 1-dimensional deformation of
$X^\infty\subseteq\P^n$ by
hypersurfaces, i.e. a family of hypersurfaces $X^\infty_s\subseteq\P^n$
of degree $d$ such that 
$X^\infty_0\simeq X^\infty$, 
$X^\infty_s$ is smooth for $s\ne 0$, and
the total space of the deformation is smooth as well.  
Consider the exact sequence of vanishing
cycles
$$
0 \to P^{n-1}(X^\infty) \to P^{n-1}(X^\infty_s)_1 \to
\oplus_j H^{n-1}(F_j)_1 \to P^n(X^\infty) \to 0\, .\leqno{(3.6)}
$$
Here $P^*(X^\infty)$ carries Deligne's $MHS,\ P^{n-1}(X^\infty_s)_1$
Schmid's MHS and\break $H^{n-1}(F_j)_1$ Steenbrink's MHS.
Moreover, the Hodge numbers of $H^{n-1}(F_j)_1$ does not depend on the
choice of the smoothing (because in a $\mu$--constant family of
hypersurface singularities, the equivariant Hodge numbers are constant), 
in particular, they are equal to the Hodge numbers provided by the 
natural local smoothing.
 The sequence
is compatible with the monodromy action, which is the identity on
$P^*(X^\infty),\ M_j$ on $H^{n-1}(F_j)$ and it will be denoted
$M_{g\ell}$ on $P^{n-1}(X_s^\infty)_1$.

\no By \cite{Sch} (resp. \cite{St.Oslo}) the weight filtration of
$P^{n-1}(X^\infty_s)_1$ (resp. $H^{n-1}(F_j)_1$) is the monodromy weight
filtration with center $n-1$ (resp.~$n$).
By the invariant cycle theorem $P^{n-1}(X^\infty)= {\rm
Ker}\,(M_{g\ell}-Id)$.  From this facts it follows that we have, for
$\ell\ge 1$
$$
\begin{array}
{l@{\hspace{2.5cm}}ll@{\hspace{3cm}}} &\dim\
Gr^W_{n-1-\ell}\ P^{n-1}(X^\infty)+ \dim\ Gr_{n-1-\ell}^W\oplus_j
H^{n-1}(F_j)_1\\[3mm] {(3.7)}\hfill&=\dim\
Gr_{n-1+\ell}^W\ \oplus_j H^{n-1}(F_j)_1-\dim\ Gr_{n-1+\ell}^W\
P^n(X^\infty).
\end{array}
$$
Now we compare
$\oplus_j H^{n-1} (F_j)_1$ and $P^n(X_t)_1$.  The
1-eigenspaces in (E.2) provide the following MHS-identification
(notice that $G^n$ has no 1--eigenvector on $P^n(X_0')$):
$$
{\rm Ker}\
(M-Id\ \vert\ P^n(X_t)_1)=\oplus^k_{j=1}\
H^{n+1}_{\{p'_j\}}\,(\XX').\leqno{(3.8)}
$$
Let $g_i$ be a local equation defining the singularity germ $(X^\infty,
p_j)$.  The smoothings of the singularities $(\XX', p'_j)$ are given
locally by the equations $$\hat g_j(x, x', x'')=g_j(x)+x'\,^2+x''\,^2,\
(x,x',x'')\in (\C^n\times\C\times\C,\ 0).$$ (cf. \cite[p. 223]{GN1},
make the change of variables $2ix'=x_0-t+x_0^{d-1},\
2x''=x_0+t-x_0^{d-1}$).
If $\hat F_j$ denotes the Milnor fiber of $\hat g_j$ then by
Sebastiani-Thom (cf.~\cite{SS}) one has $H^{n+1}(\hat
F_j)_1=H^{n-1}(F_j)_1(-1)$. \\
Let $V_1: H^{n-1}(F_j)_1\to H_c^{n-1}(F_j)_1(-1)$ be the infinitesimal
variation map of type $(-1,\,-1)$
(cf.~\cite[(2.5)]{St.Proc}), let $k:H^{n-1}_c(F_j)_1\longrightarrow
H^{n-1}(F_j)_1$ be the natural map,
$M_{j,c}:H^{n-1}_c(F_j)\longrightarrow H^{n-1}_c(F_j)$ the monodromy at
the level of cohomology with compact supports.
Define  $N_1=\log(M_1),\ N_{j,1}=\log(M_{j,1}),$ and
$N_{j,c,1}=\log(M_{j,c,1})$. Then from \cite{St.Proc} one has:
$$
H^{n-1}_{\{p_j\}}(X^\infty)={\rm Ker}\,(H^{n-1}_c(F_j)_1
\buildrel k\over\longrightarrow
H^{n-1}(F_j)_1)
$$
and similarly for $\XX'$, so one has
$$H^{n+1}_{\{p'_j\}}\ (\XX')=H^{n-1}_{\{p_j\}}\ (X^\infty)(-1).
\leqno{(3.9)}$$
Since
$N_{j,c,1}=V_1\circ k$ and $V_1$ is an isomorphism, one obtains that
$H^{n-1}_{\{p_j\}} (X^\infty)(-1)$
$={\rm Ker}\, k(-1)=
{\rm Ker}\,N_{j,c,1}(-1)\simeq {\rm Ker}\,N_{j,1}$. The last
isomorphism is given by $V_1$ and follows from the fact that
$(M_{j,c})_1\circ V_1=V_1\circ (M_j)_1$. \medskip

\no The above identities, together with (3.9) and (3.8) give:
$$
{\rm Ker}\,(N_1\ \vert\ P^n (X_t)_1)\simeq\bigoplus^k_{j=1}\
{\rm Ker}\,(N_{j,1}\ \vert\ H^{n-1}(F_j)_1).
$$
Now, (3.3) implies that
$$P^n(X_t)_1\ \buildrel n\over=
\oplus_{j=1}^k \ H^{n-1}(F_j)_1 \leqno{(3.10)}$$
(i.e. they  have the same Hodge numbers).
\medskip

\no Now consider (E.0) and (2.3). The map $N_1^0:H^n (X^0_t)_1
\longrightarrow H^n(X_t^0)_1$ is a
morphism of type $(-1,\,-1)$ and by~(2.3):
$$
{\rm Im}\ (N_1^0)={\rm Im}\ [P^n(X_t)_1\to
H^n(X_t^0)_1].\leqno{(3.11)}
$$
Via (3.4) and (3.10), identity (3.7) reads
$$(3.12)\hspace{2cm}
\begin{array}{l}
\hfill\dim\ Gr^W_{n+1+\ell}\ P^{n+1}_{X^\infty}\ (X_t)+
\dim\ Gr^W_{n+1+\ell}\ P^n(X_t)_1\\[3mm]
\quad=\dim\ Gr^W_{n+1-\ell}\ P^n(X_t)_1-
\dim\ Gr^W_{n+1-\ell}\ P^n_{X^\infty}\,(X_t).
\end{array}\hspace{1cm}
$$
Using (E.0), this equality can be rewritten as
$$
\dim\ Gr^W_{n+1+\ell}\
H^n(X_t^0)_1=\dim\ Gr^W_{n+1-\ell}\ H^n(X_t^0)_1.
$$
Together with
(3.11), this implies that the weight filtration of $H^n(X^0_t)_1$ is the
monodromy weight filtration centered at $n+1$ and its primitive spaces
can be identified with $P^{n+1}_{X^\infty}(X_t)$, i.e. for $p+q\ge n+1$
$$
p_1^{p,q}(f)=h^{p,q}(P^{n+1}_{X^\infty}\,(X_t))=
h^{n-p,n-q}(P^{n-1}(X^\infty))=
h^{n-q,n-p}\,(P^{n-1}(X^\infty)).
$$

\no {\bf Case $\xi^d=1,\ \xi\ne 1$.}
\medskip

\no From (E.0) $H^n(X^0_t)_{\xi}=P^n(X_t)_{\xi}$ and from (E.2)
$$
{\rm Ker}\ (M-\xi \,Id\ \vert\ P^n(X_t)\,)=P^n(X'_0)_{\xi}.
$$
Since the weight filtration of $P^n(X_t)_{\xi}$ is centered at $n$
(cf. \cite{Sch}),
$$
p^{p,q}_{\xi}\,(f)=h^{n-q,\ n-p}\,(P^n(X'_0)_{\xi}).
$$
\medskip

\no {\bf Case $\xi^d\neq 1$.}
\medskip

\no Given $\xi$ with $\xi^d\not=1$, from (E.0) we have
$H^n(X^0_t)_{\xi}=P^n(X_t)_{\xi}$ as $MHS$'s, and from (E.1') we have
$$
P^n(X_t)_{\xi}=(\oplus_j\ H^{n-1}(F_j)^{\oplus d-1},\ \oplus
c_{d-1}(M_j))_{\xi^{-1}} \leqno(3.13)
$$
as vector spaces.
Recall that, as vector spaces with automorphism we have
$$
\Big(H^n(F'_j), M'_j\Big)\simeq
\Big(H^{n-1}(F_j)^{\oplus d-1},\ c_{d-1}(M_j)^d\Big).\leqno(3.14)
$$
By (E.1), for $\eta\not=1$ we have an isomorphism of MHS's:
$$
\oplus_{\xi^{d}=\eta}\ P^n(X_t)_{\xi}=\oplus_j\
(H^n(F'_j),M'_j)_{\eta^{-1}}\leqno(3.15)
$$
where $(H^n(F'_j),M'_j)$ is given by the local
smoothing $\pi':(\XX', p'_j)$$\longrightarrow$$ (D', 0)$.
In some local coordinates this can be identified with
$\pi':(\YY _j,0)\longrightarrow(D',0)$, where
$$
\YY _j=\Big\{(y, y_0, t)\
\vert\ g_j(y)+t\,y_0=y_0^d\Big\}
$$
and $\pi'(y,y_0,t)=t$ (see
\cite{GN1}).  We will compute the $MHS$ of the vanishing cohomology of
$\pi'$ corresponding to eigenvalues $\neq 1$.
This computation, and the above isomorphisms will provide the result.\\
The idea of the following construction comes from a recent paper of
H. Hamm \cite{Hamm}. We thank J. Steenbrink for drawing our attention
to that paper.

\medskip \no Fix an integer $m\ge 1$ such that for any eigenvalue
$\lambda$ of $M_j$ one has $\lambda^m=1$.  Then, by (3.14), for any
eigenvalue $\lambda$ of $M'_j$ one has $\lambda^{m(d-1)}=1$.

\medskip
\no Let $(\tilde \YY _j, 0)$ be the fibre product of $\pi':(\YY _j,
0)\to(D',0)$ and $(D'',0)\to(D',0)$ where the last arrow is $s\mapsto
s^{m(d-1)}=t$ and $D''$ is again a small disk.  So, $(\tilde\YY _j, 0)$
is locally given by the equation $ g_j(y)+s^{m(d-1)}y_0=y_0^d. $
\medskip

\no Let $\pi'':(\tilde\YY _j, 0)\to(D'', 0)$
be the projection $(y_0, y, s)\mapsto s$.  Then the fiber of $\pi''$ is
the same as the fiber of $\pi'$, its monodromy is the $m(d-1)$-th power
of the monodromy of $\pi'$, in particular it is unipotent.  Let $\tilde
K$ be the real link of $(\tilde\YY _j, 0)$.
On $\tilde{K}$ we define the natural Galois action induced by
$(y_0,y,s)\mapsto (y_0,y,s\ e(\frac{1}{m(d-1)}))$.

\medskip
\no We denote by $V'_j:H^n(F'_j)\longrightarrow H^n_c(F'_j)(-1)$ the
infinitesimal variation map (\cite[(2.5)]{St.Proc}).
\no By [loc. cit., 2.6.b] one has the following exact sequence of
$MHS$'s:
$$
0\to H^n(\tilde K)\to H^n(F'_j)\buildrel V'_j\over\rightarrow
H^n_c(F'_j)(-1)\to H^{n+1}(\tilde K)\to 0\,.
$$
\noindent (3.16)\ {\bf Lemma:}\ {\em This exact
sequence is equivariant with respect to an action which at the level of
$H^n(F'_j)$ is the monodromy action $M'_j$ of $\pi'$ and on
$H^n(\tilde K)$ is the Galois action.}
\smallskip

\no {\sl Proof of the lemma:}\
If $T'_{j,geom}$ denotes the geometric monodromy of $\pi'$, then
its lifting composed with the Galois action is isotopic to the identity
(by a similar argument as in
the proof of E.1' and E.2 in \cite{GN2}). Therefore, the inverse of
$(T'_{j,geom})^*$ corresponds to the Galois action.
\smallskip

\no If $H^n(\tilde K)_\eta$ denotes the $\eta$-eigenspace $(\eta\ne 1)$
with
respect to this Galois action, then by (3.16) $H^n(\tilde K)_\eta={\rm
Ker}\,(V'_j)_{\eta}$.  Recall that $k\circ V'_j=N'_j$ where
$N'_j=\frac{1}{m(d-1)}\ \log\,(M'_j)^{m(d-1)}$ and $k$ is the natural
map $H^n_c(F'_j)\longrightarrow H^n(F'_j)$.  Since $k_{\eta}$ is an
isomorphism for $\eta\ne 1,\ {\rm Ker}\,(N'_j)_{\eta}={\rm
Ker}\,(V'_j)_{\eta}$.  Therefore, for $\eta\ne 1$,
$$
H^n(\tilde
K)_\eta={\rm Ker}\,\Big((N'_j)_{\eta}\ \vert\
H^n(F'_j)_{\eta}\Big).\leqno(3.17)
$$
The advantage of this identification
is the following:
The map $\pi'$ is defined on a hypersurface singularity germ, and there
are very few methods to compute the $MHS$ of its vanishing cohomology.
On the other hand, $\tilde K$ can be represented as the link of a
hypersurface singularity, and there are good methods to compute its
$MHS$.  Indeed, let $u:(\C^{n+2},0)\longrightarrow(\C, 0)$ be defined by
$u\,(y, y_0, s)=g_j(y)+s^{m(d-1)}\,y_0-y_0^d$.  Let $U$ be its Milnor
fiber.
Then by \cite{St.Proc} we have the exact sequence of $MHS$'s
$$
0\to
H^n(\tilde K)\to H^{n+1}_c(U)_1\buildrel k_1\over\to H^{n+1}(U)_1\to
H^{n+1}(\tilde K)\to 0\,.\leqno(3.18)
$$
Besides the monodromy action, on
this exact sequence there is a Galois action induced by the
transformation $s\mapsto s\cdot e(\frac{1}{m(d-1)})$, which commutes
with the monodromy action.  \medskip

\no Let $H^{n+1}_c(U)_{1,\eta}$ be the generalized eigenspace
corresponding to eigenvalue 1 with respect to the monodromy action and
eigenvalue $\eta$ with respect to the Galois action.  Then (3.18) gives
$$
H^n(\tilde K)_\eta={\rm Ker}\,(k_{1,\eta}\, \vert\,
H^{n+1}_c(U)_{1,\eta}).
$$
If $V_1:H^{n+1}(U)_1\longrightarrow H^{n+1}_c(U)_1(-1)$
is the infinitesimal variation map of $u$, then $V_1\circ
k=N_{c,1}$ ($N_c$ being the logarithm of monodromy acting on cohomology
with compact supports).  Since $V_1$ is an isomorphism and $N_{c,1}\circ
V_1=V_1\circ N_1$, we get\\[-15mm]
\bc
$$ \ba{l@{\hspace{1mm}}l@{\hspace{1mm}}l} {\rm Ker}\,k_{1,\eta}&=&{\rm
Ker}\,(N_{c,1,\eta}\,\vert\,H^{n+1}_c(U)_{1,\eta})\\[2mm] &=&{\rm
Ker}\,(N_{1,\eta}\,\vert\,H^{n+1}(U)_{1,\eta})(+1). \ea $$
\ec
Therefore, for $\eta\not=1$:
$$
H^n(\tilde K)_\eta={\rm Ker}\,\Big(N_{1,\eta}\,\vert\,
H^{n+1}(U)_{1,\eta}\Big)(+1).\leqno(3.19)
$$
Now, (3.17), (3.19) and (3.3) give the {\em numerical} equivalence
of $H^n(F'_j)_{\eta}$ and $H^{n+1}(U)_{1,\eta}(+1)$,  i.e.
$$
H^{pq}(F'_j)_{\eta}\buildrel n\over= H^{p+1, q+1}\ (U)_{1,\eta}.
\leqno(3.20)
$$
Let $w:(\C^2,0)\longrightarrow (\C, 0)$ be given by $w(y_0,
s)=s^{m(d-1)}y_0-y_0^d$.  Since $w$ defines an isolated quasihomogeneous
singularity, its equivariant Hodge numbers can be computed using
\cite{St.Compos} (see also \cite{AGV}).  
Then, the Sebastiani-Thom result of \cite{SS} allows
to express the equivariant Hodge numbers of $u=g_j+w$ in terms of those
of $g_j$. \medskip

\no(3.21) {\bf Proposition}:
\smallskip  {\it

\no (a) For each 1-eigenvector of $(H^{n-1}(F_j), M_j)$ of type $(p,q)$
there are $d-2$ eigenvectors of $(H^n(F'_j), M'_j)_{\ne 1}$ of type
$(p,q)$ and of eigenvalues $e(\frac{k}{d-1})$ with
$1\le k\le d-2$. \smallskip

\no (b) For each $e(i/d)$-eigenvector of type $(p,q)$ of $(H^{n-1}(F_j),
M_j)$ with $i\in\Z,\ 0<i<d$, there are $d-2$ eigenvectors
$ \{v_k\}_k$ of
$(H^n(F'_j), M'_j)_{\ne 1}$  where $k=0,\dots ,\widehat{d-
1-i},\dots d-2,\ v_k$ has eigenvalue
$e(\frac{k+i}{d-1})$ and type
$$
\Big(p+\Big[\frac{k+i}{d-1}\Big]\,,\
q+1-\Big[\frac{k+i}{d-1}\Big]\Big).
$$
(c) For each
$e(\gamma)$-eigenvector of $(H^{n-1}(F_j), M_j)$ with
$d\gamma\not\in\Z$ and  $\gamma\in (0,1)$,
and type $(p,q)$, there are $d-1$ eigenvectors $\{w_k\}_k$ of
$(H^n(F'_j),
M'_j)_{\ne 1}$  where $k=0,\dots,\,d-2\,,\ w_k$ has eigenvalue
$e(\gamma+\frac{k+\gamma}{d-1})$ and type
$$
\Big(p+\Big[\gamma+\frac{k+\gamma}{d-1}\Big],\
q+1-\Big[\gamma+\frac{k+\gamma}{d-1}\Big]\Big).
$$
(d) All eigenvectors
of $(H^{n-1}(F'_j), M_j)$ corresponding to eigenvalues $\ne 1$ are
obtained from a),~b) or~c) above. } \medskip

\no{\sl Proof:}
In order to have a simpler form of Sebastiani-Thom, we will work
with spectral pairs.  A $\lambda$-eigenvector $v$ of type $(p,q)$ of
$(H^{n-1}(F_j), M_j)$ gives a pair $S_v=(\alpha, p+q-s)$ where
$\lambda=e(-\alpha),\ n-1+[-\alpha]=p$ and $s=0$ if $\alpha\not\in\Z$
and $s=1$ if $\alpha\in\Z$  (cf. \S 6). \medskip

\no One has $Spp(w)=\sum_{(i,j)}
(\frac{i}{d}+\frac{j}{md}-1,1)$,
where the residue classes of the monomials $\{y_0^{i-1}s^{j-1}\}$
form a base in the Milnor algebra $\C[y_0,s]/(\partial w)$. A possible
choice
of the monomials is the following: $(i,j)$ either satisfies
$1\leq i\leq d$
and  $1\leq j\leq m(d-1)-1$, or $(i,j)=(1,m(d-1))$. Then $v$
contributes to
$Spp(u)$ with pairs
$$
\sum_{(i,j)}\ \Big(\alpha+\frac{i}{d}+\frac{j}{md}\,,\
p+q-s+2\Big).
$$
This gives eigenvectors of $H^{n+1}(U)_{1,\ne 1}$ if and
only if
$$
\alpha+\frac{i}{d}+\frac{j}{md}\in\Z\quad{\rm and}\quad
\frac{j}{m(d-1)}\not\in\Z.\leqno(3.22)
$$
One has to search for solutions
of (3.22) with the restrictions
$$
1\le i\le d\ ,\quad 1\le j\le
m(d-1)-1.\leqno(3.23)
$$
In case (a) (resp.  (b), resp.  (c)) one has
$S_v=(n-p-1,\ p+q-1)$ ( resp.\ $S_v=(-i/d-p+n-1,\ p+q)$,
resp. $S_v=(-\gamma-p+n-1,\ p+q)\,)$.
One gets (a), (b), (c) by solving (3.22), (3.23) in each case.
The details of the computation are left to the reader.
\bigskip

\no Now we can finish the proof of theorem 3.1 (case $\xi^d\neq 1$).
Recall the identities (3.13-14-15).
\smallskip

\no{\sl Subcase 1:} Assume $\xi^{d-1}=1$,\ i.e.
$\xi=e(-\frac{\ell}{d-1})$ with $1\le\ell\le d-2$.  Since $\xi^{-1}$ is an
eigenvalue of $c_{d-1}(M_j)$ (cf. 3.13),\
$\xi^{1-d}=1$ is an eigenvalue of~$M_j$.
By (3.21.a), each 1-eigenvector of $M_j$ of type $(p,q)$ gives $d-2$
eigenvectors of $M'_j$ of type $(p,q)$, with eigenvalues
$\eta^{-1}=e(\frac{k}{d-1})$ ($1\leq k\leq d-2$).
But (in (3.15)) only one corresponds to $\xi$ (i.e. satisfies
$e(-\frac{k}{d-1})=\xi^d$), namely that with $k=\ell$.  So, each
1-eigenvector
of $M_j$ of type $(p,q)$ gives one eigenvector of $(H^n(X_t^0),\
M^\infty_f)_{\xi}$ of type $(p,q)$.  Therefore,
$$
p^{p,q}_{\xi}(f)=\sum_j\ p^{p,q}_1(X^\infty,\,p_j).
$$

\no{\sl Subcase 2:} Assume $\xi^{(d-1)d}=1$ but $\xi^{d-1}\ne 1$.  Then
$\xi=e(-\beta)$, where $\beta=\frac{\ell}{d-1}+\frac{i}{d(d-1)}$, with
$0<i<d,\ 0\le\ell\le d-2$.
Then $\alpha=\xi^{1-d}=e(\frac{i}{d})$ is an eigenvalue of $M_j$ with
$\alpha\ne 1,\ \alpha^d=1$. Using (3.21.b) one gets eigenvalues
$\eta^{-1}=e(\frac{k+i}{d-1})$ for  $M'_j$, but only one  corresponds to
$\xi$ (because $e(-\frac{k+i}{d-1})=\xi^d$ implies $k=\ell$). Thus,
from (3.21.b), (3.13) and (3.15) one has:
$$
p^{p,q}_{\xi}(f)=\sum_j\ p_{\xi^{1-d}}^{p-[\frac{\ell+i}{d-1}],\
q-1+[\frac{\ell+i}{d-1}]}(X^\infty, p_j).
$$
\no Now, notice that $\gamma:= \{ (d-1)\beta \}=\frac{i}{d}$, hence
$[\beta+\gamma]=[\frac{\ell+i}{d-1}]$.

\smallskip

\no{\sl Subcase 3: $\xi^{d(d-1)}\ne 1$}.
\smallskip

\no Follows from (3.21.c) as in the subcases above (now $k=[(d-1)\beta]$).
\medskip

\no (3.24) {\it Remark} : The proof of theorem (3.1) in case $\xi= 1$
suggests the existence of a natural perfect pairing
$$D:P^{n-1}(X^\infty_s)_1\otimes H^n(X^0_t)_1\longrightarrow \Q(-n).$$
Our proof gives this duality at the level of Hodge numbers
and we expect that there is a geometrical construction of~$D$.
We expect that $H^n(X^0_t)_1$ is not pure if $Var$ is not an isomorphism
(cf.~(3.2.b)$\,$).  This suggests that $D$ should be constructed from
the Milnor fibration at infinity (Notice that the MHS on
$P^{n-1}(X^{\infty}_s)_1$ depends on the choice of the smoothing
$\{X_s\}_s$ of $X^{\infty}$; also the MHS of $H^n(X^0_t)_1$ depends on
the lower order terms $f_{d-1}+f_{d-2}+\cdots$ of $f$.  But maybe a good
pairing of these parameter spaces gives rise to a duality $D$).
\bigskip

\no {\bf 4. The Hodge numbers of $X^\infty$ and $X'_0$.}
\bigskip

\no Theorem (3.1) gives the $MHS$ of $\oplus_{\xi^d=1}\
H^n(X^0_t)_\xi$ in terms of the $MHS$ of $X^\infty$ and $X'_0$.  In the
sequel we would like to separate the local and the global information.
We say that an invariant is local if it depends only on the type of the
hypersurface singularities $\{(X^\infty, p_j)\}^k_{j=1}$ and some
universal constant depending on $n$ and $d$.  If an invariant depends on
the positions of the singularities of $X^\infty$, then we say that it is
global.  For example, $\dim\ P^n(X^\infty)$ and $\dim\
P^{n-1}(X^\infty)$ are global, but their difference is local. \medskip

\no (4.1) We start defining some universal constants:
If $Q\in\C\,[z_1,\dots,\,z_{n+1}]$ is an homogeneous polynomial of
degree $d$ which defines a smooth hypersurface $V=V(Q)\subseteq \P^n$,
then the Hodge numbers of $V$ can be computed as follows
(cf.~\cite{Gr}): Let $(\partial Q)$ denote the ideal in
$\C\,[z_1,\dots,\,z_{n+1}]$ generated by the partial derivatives of $Q$
and set $\FF(Q)=\C\,[z_1,\dots,\,z_{n+1}]/(\partial Q)$.  Let
$\FF(Q)_\ell$ be the subspace generated by the monomials of degree
$\ell$, and set
$$
\Omega_n=\sum_{i=1}^{n+1}\ (-1)^i\ z_i\
dz_1\wedge\dots\wedge\ \widehat{dz_i}\wedge\dots\wedge\ dz_{n+1}.
$$
If $P$ is an homogeneous polynomial of degree $\deg P=(k+1)\,d-n-1$ for
some $k\ge 0$, then $P\cdot\Omega_n\,/\,Q^{k+1}$ defines a rational
differential on $\P^n$ with poles along the hypersurface $V$.  Its
residue class ${\rm Res} (P\,\Omega_n\,/\,Q^{k+1})$ lives in the
primitive cohomology $P^{n-1}\,(V)$ and one has an isomorphism
(cf.~\cite[(4.6)]{Gr})
$$
\FF\,(Q)_{(k+1)\,d-n-1}\buildrel\sim\over\longrightarrow
P^{n-1-k,\,k}\ (V)\leqno{(4.2)}
$$
given by $P\longmapsto {\rm Res}
(P\,\Omega_n\,/\,Q^{k+1})$.

\no It is well known that the Hodge numbers of the cohomology of $V$
(equivalently, $\dim\ \FF\,(Q)_\ell$) does not depend on the choice of
$Q$ as long as $V\,(Q)\subseteq\P^n$ is smooth. We
denote
$$
j_{k,s}^{n,d}=\dim\ \FF\,(Q)_{(k+1)\,d-(n+1)-s}
$$
where $0\le s<d$ and $k\ge 0$.
The next theorem gives the Hodge numbers of $P^{n-1}(X^\infty)$ in terms
of the Hodge numbers of $P^n(X^\infty)$ modulo local contributions.
After this paper was written, we have got to know of a very recent paper
of A. Dimca (\cite{DiHa}) where the same theorem is proved, under a
slightly different form.
\medskip

\no (4.4) {\bf Theorem}: {\it
\smallskip

\no Let $X^\infty,\ \{p_1,\dots,\,p_k\}, F_1,\dots ,F_k$ be as in the
previous sections. Then:
\begin{itemize}
\item[(a)] $h^{i,\,n-k-i}\ (P^{n-1}\,(X^\infty)\,)=
\sum_j\
p_1^{k+i-1,\ n-1-i}\ (X^{\infty}, p_j)$\quad if $k\ge 3$.
\item[(b)] $h^{i,\,n-2-i}\ (P^{n-1}\,(X^\infty)\,)=
\sum_j\ p_1^{i+1,\ n-1-i}\
(X^{\infty}, p_j)-h^{i+1,\ n-1-i}\ (P^n(X^\infty)\,)$.
\item[(c)] $h^{i,\,n-1-i}\ (P^{n-1}\,(X^\infty)\,)=
j^{n,d}_{i,0}-\sum_j\ \dim\
Gr^i_F\ H^{n-1}\,(X^{\infty},p_j)\newline
-\sum_j\,\sum_{q>i}\ p_1^{q,\ n-1-i}\ (X^{\infty}, p_j)+h^{i,\,n-i}\
(P^n(X^\infty)\,)+ h^{i+1,\ n-1-i}\ (P^n\,(X^\infty)\,)$.
\end{itemize} }
\medskip

\no{\sl Proof:}
\no(a) and (b) follow from (3.6) or (3.7).
For (c), consider the following exact sequence
(similarly as in (3.6), but now for all eigenvalues):
$$
0\to P^{n-1}\,(X^\infty)\to
P^{n-1} (X^\infty_s)\to\oplus_j\ H^{n-1}\,(F_j)\to
P^n\,(X^\infty)\to 0.\leqno{(4.5)}
$$
Passing to the
limit $MHS$ does not change the dimension of the terms of the Hodge
filtration, thus $\dim\ Gr_F^i\ P^{n-1}\,(X_s^\infty)= j^{n,d}_{n-1-i,0}=
j^{n,d}_{i,0}$ (cf. 4.2).
Now take $Gr_F^i$ on the sequence above and apply (a) and (b). \bigskip

\no (4.6) {\bf Corollary}:
\smallskip {\it

\no If $p+q\ge n+3$, then $p_1^{p,q}\,(f)$ is local}. \smallskip

\no {\sl Proof:} Use the theorem
above and theorem~(3.1,a).
\medskip

\no Let $X'_0=\{f_d(z_1,\dots,\,z_{n+1})+z^d_{n+2}=0\}\subseteq\P^{n+1}$
be the $d$-th fold cyclic covering of $\P^n$ branched along $X^\infty$.
The Galois action of $\Z/d\Z$ is given by $\xi*[z_1:\dots
:z_{n+2}]=[z_1:\dots :\xi\ z_{n+2}]\ (\xi^d=1)$.  Recall that
$P^*(X'_0)_\xi$ denotes the $\xi$-eigenspace of $P^*(X'_0)$ with respect
to the Galois action.

\no The singularities ${\rm Sing}\ X'_0=\{p'_1,\dots,\,p'_k\}$ are the
$d$-th suspensions of the singularities of
$X^\infty$. If $(X^\infty, p_j)$ is locally given
in $\C^n$ by $g_j(x_2,\dots,\,x_{n+1})=0$, then $(X'_0, p'_j)$ is given
by $\tilde g_j:=g_j(x_2,\dots,\,x_{n+1})+x^d_{n+2}=0$ in $\C^{n+1}$.
Let $\tilde F_j$ be the Milnor fiber of $\tilde g_j,\ \tilde M_j$ its
monodromy acting on $H^n(\tilde F_j)$.  Then on $H^n(\tilde F_j)$ there
is a Galois action (given by $\xi *
(x_2,\dots,\,x_{n+2})=(x_2,\dots,\,\xi\,x_{n+2})$ if $\xi^d=1$) and this
action commutes with $\tilde M_j$.
$H^n(\tilde F_j)$ carries a $MHS$, we consider the space of
primitive elements given by the nilpotent part of $\tilde M_j$. Let
$P^{ab}_1\,(X'_0, p'_j)$ be the space of primitive vectors of type
$(a,b)$ corresponding to the eigenvalue~1 of $\tilde M_j$.  Finally, let
$p^{ab}_1\,(X'_0,p'_j)_\xi$ be the dimension of the $\xi$-eigenspace of
the induced Galois action on $P_1^{a,b}\,(X'_0,p'_j)$. 
(In order to eliminate any confusion, we emphasize that  below in (4.7.c)
$H^n(\tilde F_j)_\xi$ is the $\xi$--eigenspace with respect to the Galois 
action).
\medskip

\no (4.7) {\bf Theorem:}  {\it
\smallskip

\no Let $\xi=e(s/d),\ 0<s<d.$
Then with the above notations, one has
\begin{itemize}
\item[(a)] $h^{i,\,n+1-k-i}\ (P^n(X'_0)_\xi)=
\sum_j\ p_1^{k+i-1,\ n-i}\ (X'_0,p'_j)
_\xi$\qquad if $k\ge 3$
\item[(b)] $h^{i,\,n-1-i}\ (P^n(X'_0)_\xi)=
\sum_j\ p_1^{i+1,\ n-i}\ (X'_0,p'_j)
_\xi-h^{i+1,\ n-i}\ (P^{n+1}(X'_0)_\xi)$.
\item[(c)] $h^{i,\,n-i}\ (P^n(X'_0)_\xi)=
j^{n, d}_{n-i,s}-\sum_j\,\dim\,Gr_F^i\
H^n(\tilde F_j)_\xi-\sum_j\, \sum_{q>i}\ p_1^{q,\ n-i}\,
(X'_0,p_j)_\xi\break
+h^{i,\ n+1-i}\ (P^{n+1}(X'_0)_\xi)+h^{i+1,\ n-i}\
 (P^{n+1}(X'_0)_\xi).$
\end{itemize}         }
\smallskip

\no{\sl Proof:}
A deformation $(f_d)_s$ of $f_d$ induces a deformation
$(f_d)_s+z^d_{n+2}$ of $f_d+z^d_{n+2}$.  Let
$V\,(\,(f_d)_s+z^d_{n+2})=X'_s\subseteq\P^{n+1}$.  Then the exact
sequence of vanishing cycles is
$$
0\to P^n(X'_0)\to P^n(X'_s)\to \bigoplus_j\ H^n(\tilde F_j)\to
P^{n+1}(X'_0)\to 0.\leqno{
(4.8)}
$$
On this sequence there are a monodromy and a Galois action which
commute.  Now the proof of~(4.7) is the same as the proof of~(4.4) if we
replace the exact sequence~(4.5) by~(4.8) (and we take the 
$\xi$-eigenspaces, with respect to the Galois action). The remaining part is
the computation of $\dim\,Gr^i_F\ P^n(X'_s)_\xi$ (for $s\not=0$ fixed).
Take $i=n-k$.  By the discussion
in~(4.1) and with the same notations, $H^{n-k,\ k}\
(P^n(\,V(Q+z^d_{n+2})\,)_\xi)$
is given by the residues of rational differentials of type $\tilde
P\,\Omega_{n+1}\,/\,(Q+z^d_{n+2})^{k+1}$ where $\deg\,\tilde
P+(n+2)=(k+1)\,d$.  Since $\Omega_{n+1}$ is homogeneous of degree one in
$z_{n+2}$, such a form lives in the $\xi$-eigenspace with respect to the
Galois action if and only if $\tilde P=P\,(z_1,\dots
,\,z_{n+1})z_{n+2}^{s-1}$, hence $P\in\FF(Q)_{(k+1)d-(n+1)-s}$. \medskip

\no (4.9) {\bf Corollary:} {\it
\smallskip

\no If $p+q\ge n+2$, then $p^{p,q}_\xi\,(f)$ is local for any
$\xi^d=1,\ \xi\ne 1$.    } \bigskip

\no In the next section we will study the polarization properties of
$H^n\,(X^0_t)$.  This is not difficult for the generalized eigenspace
$H^n\,(X^0_t)_{\ne 1}$, but for $H^n(X_t^0)_1$ we need a suplementary
construction.
\medskip

\no Notice first that if $f=f_d+\dots$\ \ is a $(*)$-polynomial, then
$\tilde f=f\,(x_1,\dots ,\,x_{n+1})+x_{n+2}^d:\C^{n+2}\to\C$ is again a
$(*)$-polynomial.
We will compare the equivariant Hodge numbers and the polarization
properties of $f$ for eigenvalue~1 and of $\tilde f$ corresponding to
 eigenvalues $\ne 1$.
In order to do this, we need a Sebastiani-Thom-type property for the
$MHS$ at infinity of $\tilde f=f+x^d_{n+2}$ (at least for eigenvalues
$\xi$ with $\xi^d=1$). This will be established now.
Let
$X''_0=\{f_d\,(z_1,\dots,\,z_{n+1})+z_{n+2}^d+z_{n+3}^d=0\}
\subseteq\P^{n+2}$
be the $d$-fold cyclic covering of $\P^{n+1}$ branched along $X'_0$.
The Galois action on $X''_0$ is given by multiplication by $
e(1/d)$ on the last coordinate $z_{n+3}$.  Then one has
\medskip

\no (4.10) {\bf Theorem}: {\it
\smallskip

\no For $\xi=e(s/d),\ 0<s<d$ and $k\in\{n,\,n-1\}$,
$$
\begin{array}{l}
h^{p+1,\ q+1}\ (P^{k+2}\,(X''_0)_\xi)=h^{p,q}\,(P^k\,(X^\infty)\,)
\\[1.5mm]
\quad+\sum_{0<t<d\atop{t+s\ne d}}\ h^{p+[\frac{s+t}{d}],\
q+1-[\frac{s+t}{d}]}\,(P^{k+1}\,(X'_0)_{
e(\frac{s+t}{d})}).
\end{array}
$$                       }
\smallskip

\no{\sl Proof:}
Similarly as in the proof of (4.7), consider a deformation $(f_d)_s$
of $f_d$, which gives deformations $(f_d)_s+z^d_{n+2},\
(f_d)_s+z^d_{n+2}+z^d_{n+3}$.  We consider the exact sequences (4.5),
(4.8) and the analogous one for $(f_d)_s+z^d_{n+2}+z^d_{n+3}$ and we
compare $\dim\,Gr^p_F$ of the corresponding terms.

\no First we concentrate our attention on the second terms in these
exact
sequences: $P^{n-1}\,(X^\infty_s),\ P^n\,(X'_s),\ P^{n+1}\,(X''_s)$.
Since $\dim\,Gr^p_F$ remains constant when we pass to a limit $MHS$,
in order to compute it we
can replace this $MHS$ by the pure Hodge structures given by a fixed
$s_0\,(s_0\ne 0)$.  Set $Q=(f_d)_{s_0},\ X=\{Q=0\}\subseteq\P^n,\
X'=\{Q+z_{n+2}^d=0\}\subseteq\P^{n+1},\
X''=\{Q+z_{n+2}^d+z_{n+3}^d\}\subseteq\P^{n+2}$. Notice that there is
also an action of $\Z /d\Z$ on $X'$ (resp. $X''$) given by
multiplication of $z_{n+2}$ (resp. $z_{n+3}$) by $e(1/d)$.
 Again by Griffiths'
construction, $P^{n+1-(q+1),\ q+1}\ (X'')$ is freely generated by the
classes $P\,\Omega_{n+2}\ (Q+z^d_{n+2}+z^d_{n+3})^{-(q+2)}$ with
$\deg\,P+n+3=d\,(q+2),\ P\in\FF\,(Q+z^d_{n+2}+z^d_{n+3})$.  Moreover,
the $\xi$-eigenspace is given by the
classes corresponding to polynomials of the form $P=z_{n+3}^{s-1}\
\tilde P\,(z_1,\dots,\,z_{n+2})$.
Write $\tilde P=\sum_{0<t<d}\ z_{n+2}^{t-1}\
P_{t-1}(z_1,\dots,\,z_{n+1})$.  Then $\deg\
P_{t-1}+t-1+s-1+n+3=d\,(q+2)$, and
$$
P^{n+1-(q+1),\ q+1}\
(X'')_\xi=\sum_{0<t<d}\ \FF(Q)_{d\,(q+2)-n-1-t-s}.\leqno{(4.11)}
$$
If $t+s<d$ then it follows that $\deg\ P_{t-1}+t+s-1+n+2=d\,(q+2)$,
hence these polynomials $P_{t-1}$ provide rational differentials
$z_{n+2}^{t+s-1}\ P_{t-1}\ \Omega_{n+1}\ (Q+z^d_{n+2})^{-(q+2)}$
generating $P^{n-(q+1),\ q+1}\ (X')_{e(\frac{s+t}{d})}$.  If $t+s>d$,
then
$$
\deg\ P_{t-1}+t+s-1-d+n+2=d\,(q+1)
$$
hence these $P_{t-1}$'s
provide forms $z_{n+2}^{t+s-1-d}\ P_{t-1}\
\Omega_{n+1}(Q+z^d_{n+2})^{-(q+1)}$ generating
$P^{n-q,\,q}\ (X')_{e(
\frac{s+t}{d})}$.
Finally, if $t+s=d$ then $\deg\ P_{t-1}+n+1=d\,(q+1)$ hence
$P_{t-1}\ \Omega_n\ Q^{-(q+1)}$ are rational forms on $\P^n$ which
generate
$P^{n-1-q,\,q}\ (X)$ via Griffith's isomorphism.  Therefore
$$
p^{p+1,\
q+1} (X'')_{e(s/d)}=p^{p,q} (X)+\sum_{0<t<d\atop{t+s\ne d}}\
p^{p+[\frac{t+s}{d}],\ q+1-[\frac{t+s}{d}]} (X')_{e(
\frac{t+s}{d})}.\leqno{(4.12)}
$$
Passing to the limit we obtain
$$
\begin{array}{l@{\hspace{2.5cm}}ll@{\hspace{3cm}}}
{(4.13)}\hfill&\dim\ Gr^{p+1}_F\ P^{n+1}(X''_s)_\xi=\dim\ Gr_F^p\
P^{n-1}\,(X_s^\infty)&\\[2mm]
&\displaystyle\quad+\sum_{0<t<d\atop{t+s\ne d}}\ \dim\
Gr_F^{p+[\frac{t+s}{d}]}\ P^n(X'_s)_
{e(\frac{t+s}{d})}.
\end{array}
$$
Now we consider the third (the local) terms in the sequences (4.5),
(4.8) and the one for $(f_d)_s+z^d_{n+2}+z^d_{n+3}$.  There is a local
analogue of the above global argument where the rational differential
forms are replaced by forms in the Brieskorn lattice. \medskip

\no (4.14) {\bf Lemma}: {\it
Consider an isolated singularity given by a map germ $g:(\C^n,
0)\to(\C,0)$ (with
local coordinates $(x_2,\dots,\,x_{n+1})\,$), set $g'=g+x^d_{n+2},\
g''=g+x_{n+2}^d+x_{n+3}^d$.  The Galois group $\Z/d\Z$ acts on the
vanishing cohomology corresponding to $g'$ (resp.~$g''$), at the level
of coordinates this action is given by multiplication of $x_{n+2}$
(resp.~$x_{n+3}$) by $e(1/d)$.  Let $H^p_\lambda (g')_\xi$ be the
$(\lambda,\xi)$-generalized eigenspace of the $p$-th vanishing
cohomology group of $g'$ where $\xi$ (resp.~$\lambda$) is an eigenvalue
of the Galois action (resp. of the monodromy), and similarly for $g''$.
This actions are compatible with the $MHS$ on the vanishing cohomology,
let $p^{pq}_\lambda (g')_\xi$ be the dimension of the space of primitive
elements of type $(p,q)$.  Then for $\xi=e(s/d)\ (0<s<d)$ and any
$\lambda$ one has
$$
p_\lambda^{p+1,\ q+1}\,(g'')_\xi=\sum_{0<t<d\atop{t+s\ne d}}\
 p_\lambda^{p+[\frac{t+s}{d}],\
q+1-[\frac{t+s}{d}]}\ (g')_{e(\frac{t+s}{d})}+p^{pq}_
\lambda\,(g).
$$                     }

\no{\sl Proof:}
Since the information provided by the equivariant primitive
cohomology Hodge numbers is the same as the one in the set of spectral
pairs, we will give the proof at the level of spectral pairs.
If $(H, T)$ is a $MHS$ with a semisimple finite action, we refer to
\cite{St.Oslo} (5.3) or \cite{Ne1}
for the definition of the set of spectral pairs $Spp\,(H, T)$
(or see section 6).

\no By the Sebastiani-Thom formula proved in \cite{SS}
(which obviously is compatible with the Galois action), we have:
$$
Spp\
H^{n+1}_\lambda\,(g'')_\xi=\sum_{\beta+\gamma\equiv \alpha( mod \ \Z)}\
S_\beta(g)\ast Spp\ H'_{e(-\gamma)}\
(x^d_{n+2}+x^d_{n+3})_\xi\leqno{(4.15)}
$$
where $\lambda=e(-
\alpha),\ S_\beta(g)=Spp\ H^{n-1}_{e(-\beta)}\,(g)$ and
$(\alpha,\omega)\ast(\alpha',\omega')=(\alpha+\alpha'+1,\
\omega+\omega'+1)$.

\no Since $h=x^d_{n+2}+x^d_{n+3}$ is homogeneous, its vanishing
cohomology is generated (in the Brieskorn lattice)
by the forms $\omega=x_{n+2}^{u-1}\ x_{n+3}^{v-1}\
(h-1)^{[-\frac{u+v}{d}]}\
d\,x_{n+2} \wedge d\,x_{n+3}$, where $0<u,v<d$. The Galois
action on this form is $e(1/d)\ast\omega=e(v/d)
\omega$ and the spectral pair provided by $\omega$ is
$(\frac{u+v}{d}-1,\ 1)$.  Therefore
$$
Spp\ H^1_{e(-\gamma)}\
(h)_\xi=\Big\{\Big(\frac{u+s}{d}-1,\ 1\Big)\Big\}
$$
where
$\gamma=\frac{u+s}{d}-1$.  So (4.15) reads: $$Spp\ H^{n+1}_\lambda\
(g'')_\xi=\sum\ S_\beta\,(g)\ast\Big(\frac{t+s}{d}-1,\
1\Big)\leqno{(4.16)}$$ where the sum is over $0<t<d,\
\beta+\frac{t+s}{d}-\alpha\in\Z$.

\no If we denote $T(a,b)$ the transformation on $\Z[\Q\,\times\,\Z]$
given
by $T\,(a,b)\ (\alpha,\omega)=(\alpha+a,\ \omega+a+b)$, then
(again by the Sebastiani--Thom theorem), we have
$$
S_\beta(g)\ast\Big(\frac{t+s}{d}-1,1\Big)=
\extra\cases{T([\frac{t+s}{d}],\,1-[\frac{t+s}{d}])\,
Spp\, H^n_{\lambda}(g')_{e(\frac{s+t}{d})}&if\ $t+s\ne d$\cr
T\,(1,1)\,Spp\,H^{n-1}_\lambda (g)&if\ $t+s=d.$}
$$
This finishes the
proof of (4.14). \smallskip

\no We return now to the proof of~(4.10).  Consider the
middle maps in the sequences~(4.5), (4.8) and the one for
$(f_d)_s+z_{n+2}^d+z^d_{n+3}$.  Both the sources and the targets of
these maps have decompositions compatible with them.
Thus their kernels and cokernels also decompose, so for $k=n-1,\ n$
one has
$$
{(4.17)}\hspace{1cm}
\dim\ Gr^{p+1}_F\ P^{k+2}(X''_0)_\xi=
\displaystyle\sum_{0<t<d\atop s+t\ne d}\\[2mm] \dim\
Gr_F^{p+[\frac{t+s}{d}]}\ P^{k+1}\,(X'_0)_{e(\frac{s+t}{d})}$$
$$+\dim\ Gr^p_F\ P^k(X^\infty).$$
If $k=n$, then (4.10)
follows since the corresponding Hodge structures are pure.
If $k=n-1$, then apply (4.7), (4.11), (4.14) and the case $k=n$.
\bigskip

\no {\bf 5. Polarization. Connection with the real Seifert form at infinity.}
\bigskip

\no (5.1) First we recall the classification of simple
$\varepsilon$-hermitian variation
structures. The basic reference is \cite{Ne1}, see also \cite{Ne3}.

\no Let $U$ be the complexification of a finite dimensional real vector
space, set $U^*={\rm Hom}_\C (U,\C)$.  We denote
$\theta:U\longrightarrow U^{**}$ the natural isomorphism $(\,\theta(u)\
(\varphi)=\varphi(u)\,)$.  A bar over an element of $U$ denotes
complex conjugation.

\no If $\varphi\in {\rm Hom}_\C\,(U, U')$ then we denote $\bar\varphi\in
{\rm Hom}_\C\,(U, U')$ the map
$\bar\varphi\,(x)=\overline{\varphi\,(\bar x)}$.  The dual map
$\varphi^*:U'^*\longrightarrow U^*$ of $\varphi$ is defined by
$\varphi^*(\psi)=\psi\circ\varphi$. \medskip

\no (5.2) {\it Definition} (\cite[(2.1)]{Ne1}):\
An $\varepsilon$-hermitian variation structure over $\C\
(\varepsilon=\pm 1)$, abbreviated in the sequel as $\varepsilon$-$HVS$,
is a system $(U, b, h, V)$ where $U$ is the complexification of a finite
dimensional real vector space and

\no (a) $b:U\longrightarrow U^*$ is a $\C$-linear morphism with
$\overline{b^*\circ\theta}=\varepsilon\,b$

\no (b) $h$ is a $b$-orthogonal automorphism of $U$, i.e.
$\overline{h^*}\circ b\circ h=b$

\no (c) $V:U^*\longrightarrow U$ is a $\C$-linear morphism with
$\overline{\theta^{-1}\circ V^*}=-\varepsilon\,V$ and $V\circ b=h-Id$.

\no If $V$ is an isomorphism, then $(U, b, h, V)$ is called simple and
in this case $h=-\varepsilon\,V\,\overline{(\theta^{-1}\circ V^*)}^{-1}$
and $b=-V^{-1}-\varepsilon\,\overline{(\theta^{-1}\circ V^*)}^{-1}$.
\medskip

\no (5.3) {\it Example:}\
Let $f:\C^{n+1}\longrightarrow\C$ be a polynomial map which admits a
Milnor fibration at infinity (e.g.~$f$ is a $(*)$-polynomial,
cf.~(2.1~b)$\,$).  Let $F$ be the fiber of the Milnor fibration, set
$U_\R=H_n (F,\R),\ U^*_\R=H^n(F, \R)$ (which is identified with $H_n(F,
\partial F,\R)$ via the perfect pairing $H_n(F, \partial F, \R)\otimes
H^n(F,\R)\longrightarrow\R$).  Let $h_\R:U_\R\longrightarrow U_\R$ be
the real (classical)
algebraic monodromy, $b_\R: U_\R\longrightarrow U^*_\R$ the
intersection form and ${\rm Var}:U^*_\R\longrightarrow U_\R$ the
variation map.  Then \ $Var$\ is an isomorphism and
$\displaystyle(U_\R, b_\R, h_\R, {\rm Var})\otimes_\R \C$ is a
$(-1)^n-HVS$ which will be denoted $\VV (f)$.
Notice, that similarly as in the local case, $Var$ is the inverse of the
(real) Seifert form of the Milnor fibration at infinity.

\no If $f$ is a $(*)$-polynomial, then the Milnor fibration is
equivalent
to the fibration of $f$ at infinity, hence $F$ is diffeomorphic to the
generic fiber of $f$ and $(h^*_\R)^{-1}$ is exactly $T_f^\infty$.
\medskip

\no (5.4) In the sequel, if $(U, b, h, V)$ is a $\varepsilon$-$HVS$
we will assume that the eigenvalues $\lambda$ of $h$ satisfy $|\lambda
|=1$ (for the general case, see~\cite{Ne1}).
In the next examples, $\JJ_k$ denotes the $k\times k$ Jordan block with
eigenvalue~1. \begin{itemize} \item [(a)] If $b$ is an isomorphism, then
$V=(h-Id)\,b^{-1}$ is determined from $h, b$.  Up to isomorphism, there
are exactly two non-degenerate $k\times k$ matrices $b$ such that $\bar
b^*=\varepsilon\,b$ and $\JJ_k^*\ b\ \JJ_k=b$.  They are determined by
the sign of the entry $b_{k,1}$ and will be denoted $b^k_\pm$.  If we
set $\varepsilon=(-1)^n$, we can take $(b^k_{\pm})_{k,1}=\pm
i^{-n^2-k+1}$.  In particular, for $\lambda\ne 1\ (\,|\lambda|=1)$ there
are exactly two simple $(-1)^n-HVS$'s:
$$
W_\lambda^k(\pm 1)=(\C^k,\
b^k_\pm , \lambda\,\JJ_k,\ (\lambda\,\JJ_k-Id)\ (b^k_\pm)^{-1}).
$$
\item[(b)] If $\lambda=1$ and $k=1$ there are two simple $(-1)^n-HVS$'s
$$
W^1_1\ (\pm 1)=(\C, 0, Id_\C,\pm i^{n^2-1})
$$
\item[(c)] If $\lambda=1$
and $k\ge 2$, then $h=\JJ_k$ and $b$ is not an isomorphism.

\no Again, there are two simple $(-1)^n-HVS$'s
$$
W_1^k(\pm 1)=(\C^k,
\tilde
b^k_\pm, \JJ_k, \tilde V^k_\pm)
$$ and they are determined by (the sign of)
$(\tilde b^k_\pm)_{k,2}=(\,(\tilde V^k_\pm)^{-1})_{k,1}=\pm
i^{-n^2-k+2}.$
\end{itemize}

\no The sign conventions are motivated by Hodge-theoretical reasons,
they correspond to the sign conventions of the polarisations of Hodge
structures.
The structure theorem of $\varepsilon$-$HVS$ is the following:
\medskip

\no (5.5) {\bf Theorem} (\cite[(2.9)]{Ne1}):  {\it
\smallskip

\no Any simple $\varepsilon$-$HVS$ decomposes as a sum of indecomposable
ones, the decomposition is unique up to isomorphism and order of
summands.  If the eigenvalues $\lambda$ of $h$ satisfy $| \lambda |=1$,
then the indecomposable structures are $W^k_\lambda (\pm 1)$ where
$k\ge 1$. \medskip   }

\no (5.6) {\it Remark:}
\no Two real variation structures are isomorphic over $\C$ if and only
if they are isomorphic as real structures (cf. \cite[(2.10)]{Ne1}.  This
implies that if $f$ is
a polynomial which admits a Milnor fibration at infinity the information
contained in $\VV\,(f)$ (see~(5.3) above) is the same as the one
contained in the real Seifert form of the Milnor fibration at infinity.

\no If $f$ is a $(*)$-polynomial, the connection between $\VV\,(f)$ and
the Hodge theoretical invariants introduced in sections 1-4 is the
following \medskip

\no (5.7) {\bf Theorem}: {\it
\smallskip

\no Let $f$ be a $(*)$-polynomial. Then:
$$
\VV\,(f)=\bigoplus_\lambda\ \bigoplus_{2n\ge a+b\ge n+s}
\ p_\lambda^{a,b}\ (f)\cdot
W_\lambda^{r+1}\ (\,(-1)^b)
$$
where $s=0$ if $\lambda\ne 1,\ s=1$ if $\lambda=1$ and $r=a+b-n-s\ge 0$.
In particular, the equivariant Hodge numbers of $H^n(X^0_t)$ determine
completely the real Seifert form of the Milnor fibration at infinity.
\medskip               }

\no{\sl Proof:}
The proof follows closely that of Theorem~(6.1) in \cite{Ne1} (which
gives the same relation in the case of hypersurface singularities).
Notice that $\VV\,(f)$ has a direct sum decomposition
$\VV\,(f)=\oplus_\lambda\ \VV_\lambda\,(f)$ given by the
eigenvalues of the monodromy $h$.  From the same argument as in
[loc. cit.], for $\lambda\not=1$ the identity
$$
\VV_\lambda\,(f)=\bigoplus_{a,b}\
p^{a,b}_\lambda\ (f)\ W_\lambda^{r+1}\ (\,(-1)^b)
$$
follows.  More precisely, by
(E.0), $H^n(X^0_t)_\lambda=P^n(X_t)_\lambda$ for $\lambda\ne 1$, in
particular $H^n(X^0_t)_\lambda$ has a non-degenerate intersection form
and its $MHS$ is polarized by $N^0_\lambda$ (in the sense of
Cattani-Kaplan, cf.~\cite{CK}).

\no Suppose now that $\lambda=1$.  Consider the polynomial
$f+x^d_{n+2}:\C^{n+2}\longrightarrow\C$.  As remarked before, it is easy
to
see that $f+x^d_{n+2}$ is again a $(*)$-polynomial.  We want to compare
$\VV\,(f+x^d_{n+2})$ with $\VV\,(f)$ and their equivariant Hodge numbers
at infinity as well.

\no Fix $\xi=e(s/d),\ 0<s<d$.  The $\varepsilon$-HVS of the
polynomial map $x_{n+2}\longmapsto x^d_{n+2}$ is
$\oplus_{\eta^d=1\atop \eta\ne 1}\ W^1_\eta(+1)$
(with $\varepsilon=(-1)^0=+1$).
By \cite{Ne2} (which solves the Sebastiani--Thom problem for Seifert
form at infinity associated with ``good'' polynomials) one has
$$
\begin{array}{l@{\hspace{1.5cm}}ll@{\hspace{4.5cm}}}
{(5.8)}\hfill&\quad\displaystyle\VV_\xi(f+x^d)=(-1)^{n+1}
\Big[\VV(f)\otimes\Big(\bigoplus_{\eta^d=
1\atop\eta\ne 1} W^1_\eta(+1)\Big)\Big]_\xi&\\[7mm]
&\displaystyle= (-1)^{n+1}\VV_1(f)\otimes W^1_\xi(+1)+
(-1)^{n+1}\bigoplus_{\eta\ne
1\atop\eta\ne\xi}\
\VV_{\eta^{-1}\xi}(f)\otimes W^1_\eta(+1).&
\end{array}
$$
By (5.5) one can write
$$
\begin{array}{l@{\hspace{1mm}}l@{\hspace{1mm}}l}
\VV_\lambda (f)&=&\displaystyle\bigoplus_{a,b} c^{a,b}_\lambda
(f)\ W_\lambda^{a+b-n-s+1}(\,(-1)^b)\\[5mm]
\VV_\xi (f)&=&\displaystyle\bigoplus_{a,b} c^{a,b}_\xi(f+x^d)\ W_\xi^{a+b-n}(\,(-1)^b).
\end{array}
$$
By (5.8) and \cite[(3.4)]{Ne3}, after a computation one gets
$$
c^{a,b}_\xi(f+x^d)=\sum_{0<t<d}\ c^{a+[\frac{s-t}{d}],\
b+[\frac{t-s}{d}]}_{e(
\frac{s-t}{d})}\ (f).\leqno{\bf(5.9)}$$
But from the case $\lambda\ne 1$ we have
$$\begin{array}{l@{\hspace{1mm}}l@{\hspace{1mm}}l}
c^{a,b}_\xi(f+x^d)&=&p^{a,b}_\xi(f+x^d)\quad{\rm and}\\[2mm]
c^{a,b}_\lambda(f)&=&p^{a,b}_\lambda(f)\quad{\rm for\ any\ }
\lambda\ne 1.
\end{array}
$$
So (5.9) determines $c^{a,b}_1(f)$ in terms of equivariant Hodge
numbers.  This numbers can be identified with the Hodge numbers of
$X^\infty,\ X'_0,\ X''_0$ by~(3.1).
Therefore, (5.9) and (4.10) (the latter applied for $k=n-1$) give
$$
c_1^{a,b}(f)=h^{n-a,\ n-b}(\,P^{n-1}(X^\infty)\,)=p_1^{a,b}(f)
$$
(details are left to the reader).  This ends the proof.
\medskip

\no (5.10) The formula given in this theorem is the exact analogue
of the corresponding formula for isolated hypersurface singularities
(cf.~\cite{Ne1}).  Thus, roughly speaking, any connection between Hodge
numbers and topological invariants which holds for isolated hypersurface
singularities holds also for $(*)$-polynomials if we replace the $MHS$
of the local vanishing cohomology by the $MHS$ at infinity, the Milnor
fiber by the generic fiber and the local monodromy by the monodromy at
infinity.  We resume and ilustrate this in the next theorem.\medskip

\no (5.11) {\bf Theorem}:  {\it
\smallskip

\no Let $f:\C^{n+1}\longrightarrow\C$ be a $(*)$-polynomial.
\smallskip

\no a) $H^n(X^0_f)_{\ne 1}$ and $H^n_c(X^0_f)_{\ne 1}$ carry $MHS$'s
which are dual with respect to $\Q(-n)$, the weight filtration of
$H^n(X^0_f)_{\ne 1}$ is the monodromy weight filtration of $N^0_{\ne 1}$
centered at $n$.
This space has a natural $(-1)^n$-symetric polarization form (the
Poincar\'e dual of the intersection form
$\langle\,\cdot\,,\,j\,\cdot\,\rangle$) and it is polarized by the
infinitesimal monodromy $N^0_{\ne 1}$ (in the sense of Cattani-Kaplan).
\medskip

\no b) $H^n(X^0_f)_1$ and $H^n_c(X^0_f)_1$ carry $MHS$'s which are dual
with respect to $\Q\,(-n)$.  The weight filtration of $H^n(X^0_f)_1$ is
the monodromy weight filtration of~$N^0_1$
centered at $n+1$.\ The infinitesimal variation map (\cite[(2.5)]{St.Proc}
$V_1:H^n(X^0_f)_1\longrightarrow H^n_c(X^0_f)_1(-1)$ is an isomorphism
of MHS's of type $(-1,\,-1). \linebreak H^n(X^0_f)_1$ has a natural
$(-1)^{n+1}$-symmetric
polarization form given by $-\langle V_1 \,\cdot\,,\,
\cdot\,\rangle$ and it is polarized by the monodromy $N^0_1$.} \medskip

\no{\sl Proof:}
Case a) follows from (E.0) and Schmid's results \cite{Sch}.  For b),
use the polarization properties of $W_1^{r+1}(\pm\,1)$ which
behaves as in the local case.  For the local case see \cite{St.Proc}
and \cite{St.Compos} or \cite{Ne4}.\quad \medskip

\no (5.13) {\bf Corollary}:\ {\em The equivariant signature
$\sigma_{\lambda}(f)$ of the generic fiber of $f$, with respect to the
monodromy at infinity, is given by:
$$\sigma_{\lambda}(f)=
\sum_{2n\ge a+b\ge n+s} \ (-1)^b p_\lambda^{a,b}\ (f)\cdot
\Bigl(1+(-1)^{a+b-n}\Bigr)/2,$$
where $s=0$ if $\lambda\not=1$, and $s=1$ if $\lambda=1$. }
\smallskip

\no For the proof, see \cite[(6.6)]{Ne1}.
\bigskip

\no {\bf 6. The spectral pairs of a $(*)$-polynomial.}
\bigskip

\no As in the local case (cf. \cite{St.Oslo,SS}, see also
\cite[13.3.A]{AGV} and \cite{Ne1}) the equivariant Hodge
numbers $h^{p,q}_\lambda$ of the $MHS$ at infinity of a polynomial $f$
can be codified in the set of spectral pairs $Spp(f)\in\Z [\Q\times\N]$
defined by
$$
Spp(f)=\sum_{(\alpha,\omega)}\
h^{n+[-\alpha],\,\omega+s_\alpha-n-[-\alpha]}_{e(-\alpha)}\
(\alpha,\omega)
$$
where $s_\alpha=0$ if $\alpha\not\in\Z$ and
$s_\alpha=1$ if $\alpha\in\Z$.
The spectrum of $f$ is defined by
$$Sp(f)=\sum\ (\alpha)\in\Z[\Q]$$
where the sum is over the spectral pairs $(\alpha,\omega)$.

\no In the local case the spectrum has three important properties:
symmetry, Sebastiani--Thom property, and semicontinuity.  In the sequel,
we will discuss this properties briefly in our global situation (i.e.
for $(*)$-polynomials).

\medskip \no The symmetry of the weight filtration gives
$h^{pq}_{\lambda}=h^{n+s-q,n+s-p}_{\lambda}$, hence the invariance of
$Spp(f)$ with respect to the transformation $(\alpha,
n+k)\leftrightarrow (\alpha+k,\ n-k)$.  The complex conjugation (i.e.
$h^{pq}_{\lambda}=h^{qp}_{\bar{\lambda}}$) composed with the first
symmetry gives the second one:
$(\alpha,\,n+k)\leftrightarrow(n-1-\alpha,\ n-k)$.  This shows also that
the spectrum is symmetric with respect to $(n-1)/2$ (for details in the
local situation, see e.g. \cite[13.3.C]{AGV}).

\no Moreover, we have also a (global) Sebastiani-Thom theorem:
if $(\alpha,\omega),\ (\alpha',\omega')\in\Q\times\N$, define
$(\alpha,\omega)\, \ast\, (\alpha',\omega')=(\alpha+\alpha'+1,\
\omega+\omega'+1),$
and extend $\ast$ to $\Z[\Q\times\N]$ by linearity.
Set:
$$S_d=\sum_{0<s<d}\ (-s/d,\,0)\in\Z[\Q\times\N].$$
\no (6.1) {\bf Theorem}:\ {\it
If $f$ is a $(*)$-polynomial of degree $d$, then}
$$Spp(f+x^d)=Spp(f)\ast S_d.$$
{\sl Proof:}
For the spectral pairs of $f+x^d$ corresponding to eigenvalue~1 the
result is easy, it uses (3.1.a) and the fact that
$$P^n(X'_0)=\oplus_{\xi^d=1,\ \xi\ne 1}\ P^n(X'_0)_\xi.$$

\no For eigenvalues $\xi\ne 1$ with $\xi^d=1$ the result follows from
(3.1.b) and (4.10) (cf.~the proof of (5.7)$\,$).  For $\xi^d\ne 1$ one
uses (3.1.b) and the local Sebastiani-Thom theorem (\cite{SS}).  In this
later case the computation is rather long and tedious.  The details are
left to the reader. 
(For a different proof, in the context of ``cohomologically tame'' polynomials,
see \cite{NS}).  \medskip

\no The main application of (6.1) is a semicontinuity property of
$Spp(f)$. If $f$ is a $(*)$-polynomial, we will consider deformations
$f_\lambda:\C^{n+1}\to \C$ of $f$ with $\lambda\in(\C,0)$ and fixed
degree~$d$. Notice that since the $(*)$ condition defines an open dense
subset in the space of polynomials in $n+1$ variables of degree $d$,
$f_{\lambda }$ is also a $(*)$-polynomial (Actually, because of (3.2.d)
it would be enough
to consider deformations $f_{\lambda}=(f_d)_\lambda+f_{d-1}+\dots$,
where $(f_d)_\lambda$ is a deformation of $f_d$).
Since in such a deformation one has $\mu^\infty(f_{\lambda\ne 0})\ge
\mu^\infty(f)$ (where $\mu^\infty(f)$ denotes the dimension of the
middle
cohomology group of the generic fibre of $f$), one cannot expect an
upper semicontinuity
for the spectrum as in the local case, but a lower semicontinuity.
\medskip

\no (6.2) {\it Definition} (cf. \cite{St.Semi}):
Given a deformation $(f)_\lambda$ of a $(*)$-polynomial $f$ and a
subset $A\subseteq\R$, we will say that $A$ is a lower semicontinuity
domain for $(f)_\lambda$ if the function which associates to
$\lambda\in(\C,0)$ the sum of the frequencies of spectral numbers of
$f_\lambda$ in $A$ (denoted by $s_A(f_\lambda)\,$) is lower
semicontinuous, i.e. $$s_A(f_{\lambda\ne 0})\ge s_A(f).$$ \smallskip

\no (6.3)\ {\bf Theorem}: {\it
\smallskip

\no Let $f$ be a $(*)$-polynomial of degree~$d$.  Every half open
interval $(\frac{k}{d},\ \frac{k}{d}+1],\ k\in\Z$, is a lower
semicontinuity domain for deformations of $f$ with fixed degree.}
\smallskip

\no{\sl Proof:}
By (E.0) and (4.5), there is a $C\in \N$ depending only on
$n,d,p$ such that $\dim\,
Gr^p_F\,H^n(X^0_t)=C-\dim\, Gr^p_F\oplus_j H^{n-1}(F_j)$.  By the local
upper-semicontinuity result (\cite{St.Semi}) the theorem follows for the
interval $(t,t+1],\ t\in\Z$.  This, together with (6.1), gives the
result similarly as in the local case (see~[loc.cit, (2.7)]).
\smallskip

\no (6.4) {\it Remarks:} a)\ Actually, by a tedious combinatorial
argument, using the corresponding local result and (3.1), one can show
that any interval $(t,t+1],\ t\in \R$ is a lower semicontinuity
domain. We do not
emphasize this proof because a more conceptual proof is in preparation,
which works for any cohomologically tame polynomial (\cite{NS}).

\no b)\ The simmetry shows that $[t,t+1)$ is also a semicontinuity
domain.

\no c)\  As in the local case, (5.12) has the following corollary: Let
$Spp\,(f)\in\Z[\Q\times\N]$ be the set of spectral pairs of a
$(*)$-polynomial.  Let $Spp_{{\rm mod}-2} (f)$ be its projection in
$\Z\,[\Q/2\,\Z\times\N]$.  The knowledge of $Spp_{{\rm mod}-2} (f)$ is
equivalent to that of the real Seifert form of the Milnor fibration of
$f$ at infinity (cf.~\cite{Ne1}).

\bigskip

\no {\bf 7. Examples.}
\bigskip

\no (7.1) {\sl Case $n=1$}. \medskip

\no Let $f\in\C [X,Y]$ be a ($\ast$)-polynomial with highest degree form
$f_d$. Write
$$
f_d =\prod _{j=1}^m \, l^{\alpha _j}_j
$$
where the $l_j\in \C [X,Y]$ are  distinct linear forms and $\sum
\alpha _j=d$. Set $\alpha=\mbox{gcd}(\alpha_1,\dots ,\alpha_m)$.
For $x\in \R$, set
$\delta (x)=1$ if $x\in\Z$, $\delta(x)=0$ if $x\not\in \Z$.
Denote
by $\lceil x\rceil$ the smallest integer bigger or equal than $x$.
\medskip

\no The equivariant Hodge numbers of $P^2(X'_0)$ can be computed
directly, if $\xi=e(s/d)$ with $0<s<d$ then $h^{1,1}(P^2(X'_0)_{\xi})=
\delta(\frac{s\alpha}{d})$, the other numbers are zero.
Then from theorem (4.7) (using also \cite[p. 305, Th\'eor\`eme 6]{AGV}
or \cite{St.Compos} for the
local terms appearing in the formulas in (4.7)) one gets all equivariant
Hodge numbers of $P^1(X'_0)$. Now one computes the equivariant Hodge
numbers of the MHS at infinity using (3.1). One gets:
\medskip

\no a) $\ p^{1,1}_1 (f)=m-1$. \medskip

\no b) If $\xi = e(-s/d)$, $0<s<d$, then
$$p^{11}_{\xi}(f)=
-\delta(\frac{s\alpha}{d})+
\sum_j\delta(\frac{s\alpha_j}{d}).
$$
$$
p^{0,1}_{\xi}(f)=-s-1+\delta (\frac{\alpha s}{d})
+\sum_j \lceil \frac{s\alpha _j}{d} \rceil.
$$
$$
p^{1,0}_{\xi}(f)=s-1+\delta (\frac{\alpha s}{d})
-\sum_j [\frac{s\alpha_j}{d}].
$$

\no c) If $\xi^d\neq 1$, put $\xi=e(-\beta)$ with $0<\beta <1$ and set
$\gamma = \{(d-1)\beta\}$. Notice that $\beta+\gamma\not\in \Z$.

\smallskip

\no If $\xi^{d-1}=1$, then $p_{\xi}^{*,*}(f)=0$.

\no If $\xi^{d-1}\not=1$ and $\beta+\gamma<1$, then $p^{10}_{\xi}(f)=0$ and
$
p^{0,1}_{\xi}(f)= \#\{j\mid \xi^{(d-1)\alpha_j}=1\}.
$
If $\xi^{d-1}\not=1$ and $\beta+\gamma>1$, then $p^{01}_{\xi}(f)=0$ and
$
p^{1,0}_{\xi}(f)= \#\{j\mid \xi^{(d-1)\alpha_j}=1\}.
$

\no (All other primitive equivariant Hodge numbers vanish.)
\medskip

\no {\it Example:}\ If $f=x^2y^2+(x+y)^3$, then  $p^{11}_1=1,\
p^{11}_{-1}=1$
hence $h^{00}_{-1}=1$ too,\ \ $p^{01}_{e(-1/6)}=p^{10}_{e(-5/6)}=2$.
In particular, the rank of the midle (co)homology of the generic fiber
is 7, and the monodromy at infinity has a Jordan block of size two.

\bigskip

\no (7.2) {\sl Zariski's examples}.
\medskip

\no Let $f_6\in\C [X,Y,Z]$ be a homogeneous polynomial of degree $6$
which defines a sextic in $\P ^2$ with six cusps and no other
singularities. Let $f=f_6+\dots$ be any $(*)$- polynomial with highest
degree form $f_6$.

\no The Hodge numbers of $X^{\infty}=\{ f_6=0 \}\subset \P^2$
are easy to compute. Those of $X'_0$, the 6-fold cyclic covering of
$\P^2$ branched along $X^{\infty}$, depend on the position of the cusps.
The equivariant Hodge numbers of $P^3(X'_0)$ can be computed with the
aid of \cite[ch.6,(4.9), see also (3.18.ii)]{Di.book} One has two
cases: \smallskip

\no {\it Case i)} If the six cusps are on a conic then
$$
h^{2,1}_{e(1/6)}(P^3(X'_0))=1=h^{1,2}_{e(5/6)}(P^3(X'_0)).
$$

\no {\it Case ii)} If the six cusps are not on a conic then
$P^3(X'_0)=0$. \medskip

\no Now from theorem (4.7) (and using again the results in
\cite{St.Compos} or \cite[p. 305, Th\'eor\`eme 6]{AGV} to
compute the local terms appearing in (4.7)), one gets the equivariant
Hodge numbers of $P^2(X'_0)$. Finally, theorem (3.1) gives the primitive
equivariant Hodge numbers of the MHS at infinity of $f$. One has:
\medskip

\no a) Both in cases (i) and (ii),
$$
p^{1,2}_1 (f)= p^{2,1}_1 (f) = 4.
$$

\no b) If $\xi =e(s/6)$ with $0<s<6$ the numbers $p^{p,q}(f)_{\xi}$
are:
\bigskip

\begin{tabular}{|c||c|c|c|c|c|}
\hline
$(p,q)$ &   $(1,1)$   & $(2,0)$ &  $(0,2)$   & $(1,2)$  & $(2,1)$ \\
\hline
\hline
$s=5$   &   $4$ (i)   &  $0$    &  $1$ (i)   & $5$ (i)  &  $0$      \\
        &   $3$ (ii)  &         &  $0$ (ii)  & $6$ (ii) &           \\
\hline
$s=4$   &   $6$       &  $0$    &  $3$       & $0$      &  $0$      \\
\hline
$s=3$   &   $7$       &  $1$    &  $1$       & $0$      &  $0$      \\
\hline
$s=2$   &   $6$       &  $3$    &  $0$       & $0$      &  $0$      \\
\hline
$s=1$   &   $4$ (i)   &  $1$ (i)   &  $0$   & $0$      &  $5$ (i)   \\
        &   $3$ (ii)  &  $0$ (ii)  &        &          &  $6$ (ii)  \\
\hline
\end{tabular}
\bigskip

\no c) Both in cases (i) and (ii), if $\xi=e(l/30)$ with $l\in\{
1,7,11,13,17,19,23,29 \}$ then $p^{1,1}_{\xi}(f)=6$. \medskip

\no All other primitive equivariant Hodge numbers vanish.
\bigskip

\begin{footnotesize}

\no Universidad de Barcelona, Dept. de Algebra y Geometr\'{\i}a.  Gran
Via,
585, 08007 Barcelona, Spain. e-mail address: rgarcia@cerber.mat.ub.es

\no The Ohio State University, 231 West 18th Avenue, Columbus OH 43210,
USA. e-mail address: nemethi@math.ohio-state.edu

\end{footnotesize}
\end{document}